\documentclass[12pt, a4paper]{article}
\usepackage{amssymb}
\usepackage{graphicx}
\usepackage{psfrag}

\newcounter{theorem}[section]
\newtheorem{defi}{\bf Definition}[subsection]
\newtheorem{lemma}{\bf Lemma}[subsection]
\newtheorem{proposition}{\bf Proposition}[subsection]
\newtheorem{corollary}{\bf Corollary}[subsection]
\newtheorem{theorem}{\bf Theorem}[subsection]
\newtheorem{obs}{\bf Remark}[subsection]
\def\bp{\noindent{\it {\bf Proof:\ }}}
\newcommand{\ep}{\hfill$\Box$}

\def\R{I\kern -0.37 em R}
\def\N{I\kern -0.37 em N}

\newcommand{\Z}{{\bf Z}}

\begin{document}

\title{Expansive homeomorphisms of the plane.}
\author{J. Groisman}
\date{April, 2008}
\maketitle
\begin{center}
{\small IMERL, Facultad de Ingenier\'{\i}a, Universidad de la
Rep\'ublica, Montevideo Uruguay\\ (e-mail: jorgeg@fing.edu.uy)}
\end{center}

\begin{abstract}

This article tackles the problem of the classification of
expansive homeomorphisms of the plane. Necessary and sufficient
conditions for a homeomorphism to be conjugate to a linear
hyperbolic automorphism will be presented. The techniques involve
topological and metric aspects of the plane. The use of a Lyapunov
metric function which defines the same topology as the one induced
by the usual metric but that, in general, is not equivalent to it
is an example of such techniques. The discovery of a hypothesis
about the behavior of Lyapunov functions at infinity allows us to
generalize some results that are valid in the compact context.
Additional local properties allow us to obtain another
classification theorem.

\end{abstract}

\section{Introduction}
The aim of this work is to describe the set of expansive
homeomorphisms of the plane with one fixed point under certain
conditions. The original question that we asked ourselves was
whether every expansive homeomorphism of the plane was a lift of
an expansive homeomorphism on some compact surface. As it is well
known, such expansive homeomorphisms were classified by Lewowicz
in \cite{L1} and Hiraide in \cite{H}. As a matter of fact, we
began by studying whether some of the results obtained in the
previously cited article could be adapted to our new context (i.e.
without working in a compact environment but having the local
compactness of the plane). In this work I study expansive
homeomorphisms with one fixed point, singular or not, and without
stable (unstable) points. The existence of a Lyapunov function
that allows us, among other things, to generalize Lewowicz's
results on stable and unstable sets will be essential. In fact, it
will also allow us to obtain a characterization of those
homeomorphisms of the plane which are liftings of expansive
homeomorphisms on $T^{2}$. The result can be tested in any given
homeomorphism $f$ provided with a suitable Lyapunov function.
Although many of the techniques used in this work are valid for
the case where there are many singularities, we leave the study of
this situation for forthcoming papers.\\ Let $f:\R^{2} \rightarrow
\R^{2}$ be a homeomorphism of the plane that admits a Lyapunov
metric function $U$, meaning $U:\R^{2} \times \R^{2} \rightarrow
\R $ continuous and positive (i.e. it is equal to zero only on the
diagonal) and $W=\Delta (\Delta U) $ positive with $\Delta
U(x,y)=U(f(x),f(y))-U(x,y)$. We define $f$ as being $U$-expansive
if given two different points of the plane $x,y$ the following
property holds: for every $k>0$, there exists $n\in \Z $ such that
$$U(f^{n}(x),f^{n}(y))>k.$$ The main objective of this work is to
describe every expansive homeomorphism $f$ with one fixed point
where some Lyapunov metric function $U$ verifies certain
conditions concerning $f$. During this work we will require the
existence of such a Lyapunov function $U$, unlike in the compact
case where expansiveness is a necessary and sufficient condition
for the existence of a Lyapunov function (see \cite{L1} ). In the
previous reference, Lewowicz classifies expansive homeomorphisms
on compact surfaces. Our main results is (Theorem \ref{TP}): {\it
A homeomorphism $f:\R^{2} \rightarrow \R^{2}$ with a fixed point
is conjugate to a linear hyperbolic automorphism if and only if it
admits a Lyapunov metric function that satisfies condition {\bf
HP} and it has not singular points. Condition {\bf HP} establishes
that: given any compact set $C$ of $\R^{2}$, the following
property holds $$ \lim_{ x\rightarrow \infty }
\frac{|V(x,y)-V(x,z)|}{W(x,y)} =0,$$ uniformly with $y,z$ in $C$
and $V=\Delta U, W=\Delta V$.}\\ Without condition {\bf HP} and
demanding other kind of conditions for $U$, different behaviors
appear. These are described in Theorem \ref{LOC}: {\it Let $f$ be
a homeomorphism of the plane with a fixed point. $f$ admits a
Lyapunov function $U$ that verifies hypothesis {\bf HL} if and
only if $f$ restricted to each quadrant determined by the stable
and unstable curves of the fixed point is conjugated (such
conjugations must preserve stable and unstable curves) either to a
linear hyperbolic automorphism or to a restriction of a linear
hyperbolic automorphism to certain invariant region.} The most
important part of condition {\bf HL} establishes that:
\begin{itemize}
\item the first difference $V=\Delta U$ verifies the following property:
given $\epsilon >0$ there exists $\delta >0$ such that if
$U(z,y)<\delta $ then $|V(x,z)-V(x,y)|<\epsilon $, $\forall x\in
\R^{2}$;

\item the second difference $W=\Delta^{2} U$ verifies the following property:
given $\delta >0$, there exists $a(\delta )>0$ such that
$W(x,y)>a(\delta )>0$ for every $x,y$ on the plane with
$U(x,y)>\delta$.

\end{itemize}
The difference between the two cases (presented in Theorem
\ref{LOC}) consists of the existence of stable and unstable curves
that do not intersect each other. In \ref{EJ} we will show
examples about the case where there are stable and unstable curves
that do not intersect each other. \\We also conjecture that if
$f:\R^{2} \rightarrow \R^{2}$ is a preserving-orientation and
fixed point free homeomorphism that admits a Lyapunov function
$U:\R^{2} \times \R^{2} \rightarrow \R$ satisfying condition {\bf
HP}, then it must be topologically conjugate to a translation of
the plane. We believe that the proof of this assertion is a
consequence of Brouwer's translation theorem (see \cite{B},
\cite{F}) and of some techniques used in this article. We leave
the study of this situation for forthcoming works. \\Regarding the
structure of the paper, we begin section \ref{PRE} by studying
some properties that are verified by a Lyapunov function
associated to a lift of an expansive homeomorphism in the compact
case, as well as to homeomorphisms conjugated to it. In section
\ref{CEI} we describe stable and unstable sets, by adapting
Lewowicz's arguments used on \cite{L1}. In Section \ref{SPG} we
show our main result. In Section \ref{OBS} we study an other
context and some examples.

\section{Preliminaries}
\label{PRE} During the course of this work we will consider
homeomorphisms of the plane which admit a Lyapunov metric function
$U$ with certain characteristics. These properties are natural
since they are verified by a Lyapunov function of a lift of an
expansive homeomorphism in the compact case. In this section we
will verify some of these properties. In \cite{Fa}, \cite{T} is
proved that every lift $F$ of an expansive homeomorphism on
compact surfaces satisfies the existence of pseudo-metrics
$D_{s}$, $D_{u}$ and $\lambda >1$ such that for every $\xi ,\ \eta
\in \R^{2}$ the following holds:
$$D_{s}(F^{-1}(\xi ),F^{-1}(\eta ))=\lambda D_{s}(\xi ,\ \eta );$$
$$D_{u}(F(\xi ),F(\eta ))=\lambda D_{u}(\xi ,\ \eta ),$$
where $D=D_{s}+D_{u}$ is a Lyapunov metric in $\R^{2} $ for $F$.\\
We will test the following properties:
\begin{itemize}
\item[(I)] {\bf Signs for $\Delta (D)$.} We shall use the notation
$B_{k}(x)$ for the connected component of set $B_{k}(x)=\{ y\in
\R^{2} \ / D(x,y)\leq k\} $, which contains $x$. For every point
$x\in \R^{2} $ and for every $k>0$, there are points $y$ in the
border of $B_{k}(x)$ such that $\Delta
D(x,y)=D(f(x),f(y))-D(x,y)>0$ and points $z$ in the border of
$B_{k}(x)$ such that $\Delta D(x,z)=D(f(x),f(z))-D(x,z)<0$. This
property will be essential to describe stable and unstable sets.

\item[(II)] {\bf Property HP.} Let $V=\Delta D$ and $W=\Delta^{2}
D$. Given any compact set $C\subset \R^{2} $, the following
property holds
$$ \lim _{\| x\| \rightarrow \infty } \frac{|V(x,y)-V(x,z)|}{W(x,y)}
=0,$$ uniformly with $y,z$ in $C$. This property will be essential
to prove the main result on this work.

\end{itemize}

\subsection{Lifted case.}
Let $D=D_{s}+D_{u}$ be the Lyapunov metric function that we
introduced at the beginning of this section.

\begin{itemize}
\item[(I)] {\bf Signs for $\Delta D$.} \bp
$$ \Delta D(x,y)=D(f(x),f(y))-D(x,y) = $$
$$ D_{s}(f(x),f(y))-D_{s}(x,y)+ D_{u}(f(x),f(y))-D_{u}(x,y)=$$
$$ (\lambda -1)D_{u}(x,y) -(1-1/\lambda )D_{s}(x,y).$$
For every point $x\in \R^{2} $ and for every $k>0$, there are
points $y$ in the border of $B_{k}(x)$ such that $D_{u}(x,y)=0$
(this is true because the stable set separates the plane).
Therefore, $\Delta D(x,y)<0$ as we wanted. A similar argument lets
us find points $z\in \R^{2}$ such that $\Delta
D(x,z)>0$.\ep

\item[(II)] {\bf Property HP.} \bp Since $$ \Delta^{2}
D(x,y)= \Delta D(f(x),f(y))-\Delta D(x,y)=$$ $$ (\lambda
-1)^{2}D_{u}(x,y) + (1-1/\lambda )^{2}D_{s}(x,y),$$ we can
conclude that $\Delta^{2} D(x,y)$ tends to infinity when $x$ tends
to infinity. Now,
$$|\Delta D(x,y)-\Delta D(x,z)|\leq $$ $$(\lambda
-1)|D_{u}(x,y)-D_{u}(x,z)|+(1-1/\lambda
)|D_{s}(x,y)-D_{s}(x,z)|\leq $$ $$(\lambda -1)D_{u}(z,y)
+(1-1/\lambda )D_{s}(z,y).$$ Then $|\Delta D(x,y)-\Delta D(x,z)|$
is uniformly bounded when points $y$ and $z$ lie on a compact set.
Then property {\bf HP} holds.\ep

\end{itemize}

\subsection{Lifted conjugated case.}
Now, let us start with the case where $f$ is conjugated to a lift
$F$ of an expansive homeomorphism on a compact surface. Let us
define a Lyapunov function for $f$ such as
$$L(p_{1},p_{2})=D(H(p_{1}),H(p_{2})),$$ where $D$ is the previous
defined Lyapunov metric function for $F$ and $H$ is a
homeomorphism from $\R^{2}$ over $\R^{2}$. It follows easily that
$L$ is a Lyapunov function for $f$ and a metric in $\R^{2}$.

\begin{itemize}

\item[(I)] {\bf Signs for $\Delta (L)$.} \bp It is clear
since
$$\Delta L(p_{1},p_{2})=\Delta D(H(p_{1}),H(p_{2})),$$
and $H$ is continuous at infinity.\ep

\item[(II)] {\bf Property HP.} \bp $$|\Delta (L)(p,q)-
\Delta (L)(p,r)|=$$ $$|\Delta D(H(p),H(q))- \Delta
D(H(p),H(r))|\leq $$
$$(\lambda -1)D_{u}(H(q),H(r)) +(1-1/\lambda )D_{s}(H(q),H(r))\leq
K,$$ since $q$ and $r$ are in a compact set and $H$ is a
homeomorphism. Since $$\Delta^{2} (L)(p,q)=\Delta^{2}
D(H(p),H(q))$$ and $H$ is continuous at infinity we conclude that
$\Delta^{2} (L)(p,q)$ tends to infinite when $p$ tends to
infinity. Then, if $f$ is conjugated to a lift of an expansive
homeomorphism on a compact surface, it admits a Lyapunov function
$L$ such that condition {\bf (HP)} holds.\ep
\end{itemize}
\section{Stable and unstable sets.}
\label{CEI}

In this section we will stay close to the arguments used by
Lewowicz in \cite{L1}, Sambarino in \cite{SAM} and Groisman in
\cite{JG}. We have to adapt them for our non-compact context. We
will work with the topology induced by a Lyapunov function $U$ and
define the $k$-stable set in the following way:
$$S_{k}(x)=\{y\in \R^{2}:\ U(f^{n}(x),f^{n}(y))\leq k, \forall n\in \N \}.$$
Similar definition for the $k$-unstable set $U_{k}$. Let $f$ be a
homeomorphism of the plane that admits a Lyapunov function
$U:\R^{2}\times\R^{2}\rightarrow\R$ such that the following
properties hold:
\begin{itemize}

\item[(1)] {\bf $U$ is a metric in $\R^{2} $ and induces the same
topology in the plane as the usual metric.} Observe that given any
Lyapunov function it is possible to obtain another Lyapunov
function that verifies all the properties of a metric except,
perhaps, for the triangular property.

\item[(2)] {\bf Existence of both signs for the first difference
of $U$.} For each point $x\in \R^{2} $ and for each $k>0$ there
exist points $y$ and $z$ on the border of $B_{k}(x)$ such that
$V(x,y)=U(f(x),f(y))-U(x,y)>0$ and $V(x,z)=U(f(x),f(z))-U(x,z)<0$,
respectively.

\end{itemize}

\begin{obs}
A homeomorphism $f$ that admits a Lyapunov function $U$ defined at
$\R^{2}\times\R^{2}$ is U-expansive. This means that given two
different points of the plane $x,y$ and given any $k>0$, there
exists $n\in \Z $ such that $$U(f^{n}(x),f^{n}(y))>k.$$
\end{obs}

\bp
Let $x$ and $y$ be two different points of the plane such that
$V(x,y)>0$. Since $\Delta V>0$, then $V(f^{n}(x),f^{n}(y))>V(x,y)$
holds for $n>0$. This means that $U(f^{n}(x),f^{n}(y))$ grows to
infinity, since
$$U(f^{n}(x),f^{n}(y))=U(x,y)+\sum_{j=1}^{n}
V(f^{j}(x),f^{j}(y))>$$
$$U(x,y)+nV(x,y).$$ Thus, given $k>0$ there exists $n\in \N $ such that
$$U(f^{n}(x),f^{n}(y))>k.$$ By using similar arguments we can prove the case
when $V(x,y)=U(f(x),f(y))-U(x,y)<0$. If $V(x,y)=0$, then
$V(f(x),f(y))>0$ and this is precisely our first case.
\ep

\begin{defi}
Let $f:\R^{2} \rightarrow \R^{2}$ be a homeomorphism of the plane
that admits a Lyapunov metric function $U$. A point $x\in \R^{2} $
is a stable (unstable) point if given any $k'>0$ there exists
$k>0$ such that for every $y\in B_{k}(x)$, it follows that
$U(f^{n}(x),f^{n}(y))<k'$ for each $n\geq 0$ ($n\leq 0$).
\end{defi}

\begin{obs}
Property $(2)$ for $U$ implies the non-existence of stable
(unstable) points.
\end{obs}

\bp Given the existence of both signs for
$V(x,y)=U(f(x),f(y))-U(x,y)$ in any neighborhood of $x$, we can
state that for each $k>0$, there exists a point $y$ in $B_{k}(x)$
such that $V(x,y)>0$. Since $\Delta V>0$, we can state that
$V(f^{n}(x),f^{n}(y))>V(x,y)$ for $n>0$, so $U(f^{n}(x),f^{n}(y))$
grows to infinity. Thus, there are no stable points. We can use
similar arguments for the unstable case.\ep

\begin{obs}
There does not exist $x,y\in \R^{2}$ and $n\in \Z $ such that
$$U(f^{n+1}(y),f^{n+1}(x))>U(f^{n}(y),f^{n}(x))$$ and
$$U(f^{n+1}(y),f^{n+1}(x))>U(f^{n+2}(y),f^{n+2}(x)).$$
\end{obs}

\bp Suppose that there exist two different points $x,y\in \R^{2}$
and $n\in \N $ such that they do not verify the thesis. Since $$
\triangle (\triangle U)(f^{n}(x),f^{n}(y))=$$
$$U(f^{n+2}(x),f^{n+2}(y))-2U(f^{n+1}(x),f^{n+1}(y))+U(f^{n}(x),f^{n}(y)),$$
we have that $$ \triangle (\triangle U)(f^{n}(x),f^{n}(y))<
U(f^{n+1}(x),f^{n+1}(y))-2U(f^{n+1}(x),f^{n+1}(y))+$$
$$U(f^{n+1}(x),f^{n+1}(y))=0,$$ which is not possible.\ep

\begin{lemma}
Let $A$ be an open set of $\R^{2} $ with $x\in A\subset B_{k}(x)$.
There exists a compact connected set $C$ with $x\in C\subset
\overline{A} $, $C\cap \partial (A)\neq \emptyset $ such that, for
all $y\in C$ and $n\geq 0$, $U(f^{n}(x),f^{n}(y))\leq k $ holds.
\end{lemma}

\bp Suppose that there exists $N>0$ such that for each compact
connected set $D\subset \overline{A}$ that joins $x$ with the
border of $A$, there exists $z\in D$ and $n$ with $0\leq n\leq N$
such that $U(f^{n}(x),f^{n}(z))>k $. Otherwise, for each $n\geq 0$
we would find $D_{n}\subset \overline{A}$ that joins $x$ with the
border of $A$ such that for every $y\in D_{n}$,
$U(f^{m}(x),f^{m}(y))\leq k $, $0\leq m\leq n$. Then $$D_{\infty
}=\bigcap_{n=0}^{\infty } \overline{\left( \bigcup_{j=n}^{\infty
}D_{j}\right) } ,$$ is a connected compact set that satisfies our
assertion. Let us go back to the prior assumption. Consider a
point $y$ in the border of $B_{k}(f^{n}(x))$ that belongs to the
region where $V(f^{n}(x),y)=U(f^{n+1}(x),f(y))-U(f^{n}(x),y)<0$,
this means that $U(f^{n+1}(x),f(y))<k$. Then
$U(f^{n-1}(x),f^{-1}(y))>k$, because otherwise, we would
contradict the previous remark. So, $\overline{B_{k}(f^{n}(x))}$
contains points $y$ such that $f^{-1}(y)$ does not belong to
$\overline{B_{k}(f^{n-1}(x))} $. Let us take $n>N$ and a point
$y$, like we did before. Let us base our reasoning in the
connected component of $B_{k}(f^{n}(x))$ which contains
$f^{n}(x)$. Let $a:[0,1]\rightarrow B_{k}(f^{n}(x))$ be an arc
such that $a(0)=f^{n}(x)$, $a(1)=y$, and let $s^{*}$ be the
supremum of $s\in [0,1]$ such that, for all $u\in [0,s]$,
$f^{p-n}(a(u))\in B_{k}(f^{p}(x))$, for all $0\leq p\leq n$ and
$f^{-n}(a(u))\in A$. Since $\Delta^{2} (U)=W>0$,
$f^{p-n}(a(s^{*}))\in B_{k}(f^{p}(x))$ for all $0\leq p\leq n$ and
then $f^{-n}(a(s^{*}))$ belongs to the border of $A$. Hence
$f^{-n}(a([0,s^{*}]))$ is a connected compact set that joins $x$
with the $\partial A$ and remains inside $B_{k}(f^{n}(x))$ for all
$0\leq n \leq N$, which contradicts our assumption.\ep For
$x\in \R^{2}$ and $k>0$, let $S_{k}(x)$ be the $k -$stable set for
$x$, defined by $$S_{k}(x)=\{ y\in \R^{2} :\
U(f^{n}(x),f^{n}(y))\leq k ,\ n\geq 0\} .$$

\begin{lemma}
Let us consider $0<k'<k $. There exists $\sigma >0$ such that if
$y\in S_{k}(x)$ and $U(x,y)<\sigma $, then $y\in S_{k'}(x)$.
\end{lemma}

\bp Let us consider $\sigma \leq k'$ and $y$ such that
$U(x,y)<\sigma $ and $y\in S_{k}(x)$. Then we state that
$V(f^{n}(x),f^{n}(y))<0$ for each $n>0$. If there exists $n_{0}>0$
such that $V(f^{n_{0}}(x),f^{n_{0}}(y))>0$, then
$V(f^{n}(x),f^{n}(y))>0$ for each $n>n_{0}$ since $W=\Delta
(V)>0$. But then, $U(f^{n}(x),f^{n}(y))$ is unbounded, which
contradicts the fact that $y\in S_{k}(x)$. But if
$V(f^{n}(x),f^{n}(y))<0$ for each $n>0$, we have that
$U(f^{n}(x),f^{n}(y))<U(x,y)<\sigma \leq k'$ for all $n>0$, and
then this proves the lemma.\ep Let $C_{k}(x)$ ($D_{k}(x)$)
be the connected component of $S_{k}(x)$ ($U_{k}(x)$) that
contains $x$.
\begin{lemma}
$C_{k}(x)$ is locally connected at $x$.
\end{lemma}
\bp (See Lemma $2.3$, \cite{L1})\ep
\begin{corollary}
For each $x$ in $\R^{2} $, $C_{k }(x)$ is connected and locally connected.
\end{corollary}
\bp (See \cite{L1})\ep
\begin{corollary}
For any $x$ in $\R^{2} $ and any pair of points $p$ and $q$ in
$C_{k}(x)$, there exists an arc included in $C_{k}(x)$ that joins
$p$ and $q$.
\end{corollary}
\bp (See Topology, Kuratowski, \cite{K}, section $50$)
\ep

\begin{defi}
We say that $p\in \R^{2}$ has local product structure if a map
$h:\R^{2} \rightarrow \R^{2}$ which is a homeomorphism over its
image ($p\in Im(h)$) exists and there exists $k>0$ such that for
all $(x,y)\in \R^{2}$ it is verified that $h(\{x\} \times
\R)=C_{k}(h(x,y))\cap Im(h)$ and $h(\R \times
\{y\})=D_{k}(h(x,y))\cap Im(h)$.
\end{defi}

\begin{figure}[htb]
\psfrag{alpha}{$\alpha $} \psfrag{Wualpha}{$W^{u}(\alpha )$}
\psfrag{Wsalpha}{$W^{s}(\alpha )$} \psfrag{phitn}{$\varphi (t_{n})
$} \psfrag{Wsx}{$W^{s}(x)$}
\begin{center}
\includegraphics[scale=0.25]{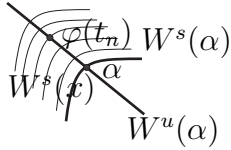}
\caption{\label{fig24}Separates the plane}
\end{center}
\end{figure}

\begin{proposition}
Except for a discrete set of points, that we shall call singular,
every $x$ in $\R^{2} $ has local product structure. The stable
(unstable) sets of a singular point $y$ consists of the union of
$r$ arcs, with $r\geq 3$ that meet only at $y$. The stable
(unstable) arcs separate unstable (stable) sectors.
\end{proposition}
\bp (See Section $3$, \cite{L1}) \ep
\begin{obs}
The neighborhood's size where there exists a local product
structure may become arbitrarily small. However we are able to
extend these stable and unstable arcs getting curves that we will
denote as $W^{s}(x)$ and $W^{u}(x)$, respectively. If two points
$y$ and $z$ belong to $W^{s}(x)$ ($W^{u}(x)$), then
$U(f^{n}(y),f^{n}(z)<K$ for some $K>0$ and for all $n\geq 0$
($n\leq 0$).
\end{obs}
The following lemmas refer to these stable and unstable curves.
\begin{lemma}
Let $f$ be a homeomorphism of the plane which verifies the
conditions of this section. Then stable and unstable curves
intersect each other at most once.
\end{lemma}
\bp If they intersect each other more than once, we would
contradict expansiveness: if two different points $x$ and $y$
belong to the intersection of a stable and an unstable curve, then
there exists $k_{0}>0$ such that $U(f^{n}(x),f^{n}(y)<k_{0}$ for
all $n\in \Z $.\ep
\begin{lemma}
Every stable (unstable) curve separates the plane.
\end{lemma}
\bp Let $\varphi :(-\infty ,+\infty )\rightarrow \R^{2} $ be a
parametrization of a stable curve $W^{s}(x)$, such that $\varphi
(0)=x$. We will first prove that $\lim _{t\rightarrow \pm \infty
}\varphi (t)=\infty $, i.e. given any closed neighborhood $B$ of
$x$ there exists $t^{*}>0$ such that $\varphi (t)$ does not belong
to $B$ for each $t\in (\R -[-t^{*},t^{*}])$. We will work with
$t>0$. Arguments for $t<0$ are identical. Let us assume that there
exists $B$ such that for each $n\in \N $ there exists $t_{n}>n$
with $\varphi (t_{n})\ \in \ B$. Let us consider the set $\{
\varphi (t_{n}):\ n\in \N \}$. This set accumulates at a point
$\alpha $ of $\overline{B}$. Let us take a large enough $p\in \N$
such that $W^{s}(\varphi (t_{p}))$ is close enough to
$W^{s}(\alpha )$ and such that $W^{u}(\alpha )$ intersects
$W^{s}(\varphi (t_{p}))$. But then $W^{u}(\alpha )$ intersects
$W^{s}(\varphi (t_{n}))$ for each $n\geq p$, which contradicts the
previous lemma (see figure \ref{fig24}).
Then we proved that $\lim _{t\rightarrow \pm \infty }\varphi
(t)=\infty $. This implies that $\R^{2}-W^{s}(x)$ has more than
one connected component. Since a stable curve can not auto
intersect (because there are no stable points), then there exist
exactly two components in the complement of $W^{s}(x)$.\ep

\section{Main section.}
\label{SPG}

In this section, we will prove the main result of this work. Let
$f$ be a homeomorphism of the plane which verifies the following
conditions:

\begin{itemize}

\item $f$ has a fixed point;

\item the quadrants determined by the stable and the unstable
curves of the fixed point are $f$-invariant;

\item $f$ admits a Lyapunov metric function $U:\R^{2} \times
\R^{2} \rightarrow \R $. This metric induces in the plane the same
topology as the usual distance;

\item $f$ has no singularities;

\item for each point $x\in \R^{2} $ and any $k>0$, there exist
points $y$ and $z$ in the border of $B_{k}(x)$ such that
$V(x,y)=U(f(x),f(y))-U(x,y)>0$ and $V(x,z)=U(f(x),f(z))-U(x,z)<0$,
where $B_{k}(x)=\{ y\in \R^{2} \ / U(x,y)\leq k\} $.

\item {\bf Property HP.} Given any compact set $C\subset \R^{2} $ the
following property holds: $$ \lim _{\| x\| \rightarrow \infty }
\frac{|V(x,y)-V(x,z)|}{W(x,y)} =0,$$ uniformly with $y,z$ in $C$.

\end{itemize}

\subsection{Previous lemmas.}

\begin{lemma}
\label{TG}

Let $f:\R^{2} \rightarrow \R^{2}$ be a homeomorphism of the plane
under the hypothesis described at the beginning of this section.
If the unstable (stable) curve of a point $x$ intersects the
stable (unstable) curve of the fixed point, then the stable
(unstable) curve of $x$ intersects the unstable (stable) curve of
the fixed point.
\end{lemma}

\begin{figure}[htb]
\psfrag{p}{$p$}  \psfrag{Wsy}{$W^{s}(y)$}
\psfrag{Wsfy}{$W^{s}(f(y))$} \psfrag{f(y)}{$f(y)$} \psfrag{y}{$y$}
\psfrag{x}{$x$} \psfrag{B}{$B$}\psfrag{f-1B}{$f^{-1}(B)$}
\psfrag{(I)}{$(I)$} \psfrag{(II)}{$(II)$} \psfrag{Wsp}{$W^{s}(p)$}
\psfrag{Wux}{$W^{u}(x)$}
\begin{center}
\includegraphics[scale=0.25]{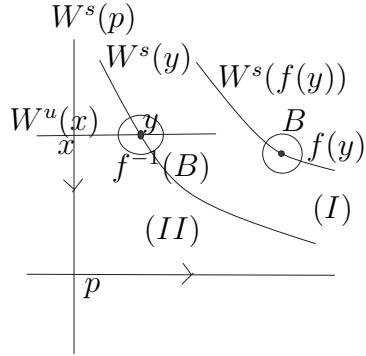}
\caption{\label{fig131}Invariant stable I}
\end{center}
\end{figure}

\bp Let us assume that there exists a point $x$ in the stable
curve ($W^{s}(p)$) of the fixed point, such that there exists a
point $y\in W^{u}(x)$ that verifies $W^{s}(y)\cap
W^{u}(p)=\emptyset $. As there are no singular points, let $y$ be
the first point in $W^{u}(x)$ such that its stable curve does not
intersect the unstable curve of the fixed point. (See figure
\ref{fig131} )

We will prove that $W^{s}(y)$ is $f$-invariant. Let us concentrate
on one of the quadrants determined by the stable and the unstable
curves of the fixed point. Since $W^{s}(y)$ separates the first
quadrant, denote by zone $(I)$ the region that does not include
the fixed point $p$ in its border, and zone $(II)$ the other
region. We will divide the proof in three steps:
\begin{itemize}
\item {\bf $f(y)$ belongs to zone $(I)$.}\\ If $f(y)$ is included
in zone $(I)$ (see figure \ref{fig131}), then $W^{s}(f(y))\cap
W^{u}(p)=\emptyset $ and the same happens to every point in a
neighborhood $B$ of $f(y)$. But then $f^{-1}(B)$ is an open set
that contains $y$ and has the property that the stable curves of
their points do not cut the unstable curve of the fixed point
(since if some point $x$ of $f^{-1}(B)$ verified that its stable
curve cuts the unstable curve of the fixed point, then $f(x)\in B$
would have the same property, which is a contradiction). Then $y$
would not be the first point of $W^{u}(x)$ such that its stable
curve does not cut the unstable curve of the fixed point (see
figure \ref{fig131}). \item {\bf $f(y)$ belongs to zone $(II)$ and
$W^{u}(f(y))\cap W^{s}(y)= \{ q\} $.}\\ Consider the simple closed
curve $J$ determined by the unstable arc $xy$, the stable arc
$yq$, the unstable arc $f(x)q$ and the stable arc $xf(x)$ (see
figure \ref{fig221}).

\begin{figure}[htb]
\psfrag{p}{$p$} \psfrag{Wup}{$W^{u}(p)$} \psfrag{Wsp}{$W^{s}(p)$}
\psfrag{Wuy}{$W^{u}(y)$} \psfrag{Wsy}{$W^{s}(y)$}
\psfrag{Wufy}{$W^{u}(f(y))$} \psfrag{f(y)}{$f(y)$} \psfrag{y}{$y$}
\psfrag{x}{$x$} \psfrag{q}{$q$} \psfrag{f(x)}{$f(x)$}
\begin{center}
\includegraphics[scale=0.25]{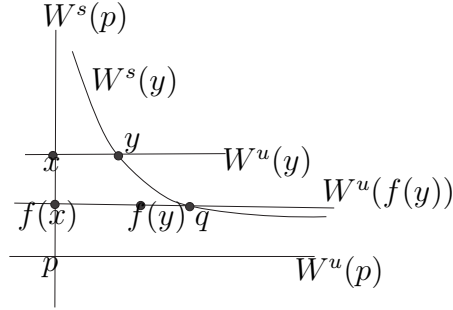}
\caption{\label{fig221}Invariant stable II}
\end{center}
\end{figure}

$W^{s}(f(y))$ intersects the bounded component determined by $J$
and intersects $J$ in $f(y)$. We will prove that $W^{s}(f(y))$
does not intersect $J$ in any other point, which be a
contradiction since $W^{s}(f(y))$ separates the plane. Indeed,
$W^{s}(f(y))$ can not cut the stable arcs of $J$ because they are
different stable curves. $W^{s}(f(y))$ can not cut the unstable
arc $f(x)q$ in a point different from $f(y)$ because we would have
a stable curve and an unstable curve cutting each other in two
different points, and we already have proved this is not possible.
If the intersection point was $f(y)$, we would have a closed curve
of a stable arc which would imply the existence of a stable point,
and this is not possible. Finally, the intersection of
$W^{s}(f(y))$ with $J$ can not belong to the unstable arc $xy$
since in this case $y$ would not be the first point of that arc
such that $W^{s}(y)$ does not intersect the unstable curve of the
fixed point.\item {\bf $f(y)$ belongs to zone $(II)$ and
$W^{u}(f(y))\cap W^{s}(y)= \emptyset $.}\\ Let us consider a
sequence $(y_{n})$ in $W^{u}(y)$ converging to $y$ such that the
stable curve of each point $y_{n}$ cuts the unstable curve of the
fixed point.
\begin{figure}[htb]
\psfrag{p}{$p$} \psfrag{Wup}{$W^{u}(p)$} \psfrag{Wsp}{$W^{s}(p)$}
\psfrag{Wuy}{$W^{u}(y)$} \psfrag{Wsy}{$W^{s}(y)$}
\psfrag{Wsyn}{$W^{s}(y_{n})$} \psfrag{Wufy}{$W^{u}(f(y))$}
\psfrag{f(y)}{$f(y)$} \psfrag{y}{$y$} \psfrag{x}{$x$}
\psfrag{yn}{$y_{n}$}\psfrag{zn}{$z_{n}$}
\begin{center}
\includegraphics[scale=0.25]{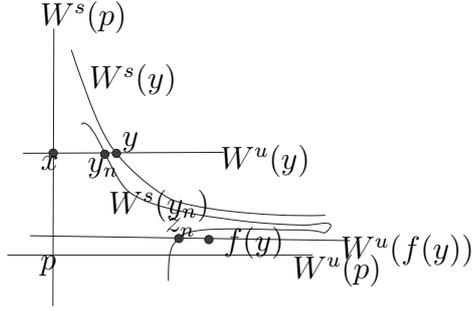}
\caption{\label{fig121}Invariant stable III}
\end{center}
\end{figure}
The behavior of the stable curves of points $y_{n}$ must be the
one shown in the figure \ref{fig121}: as $n$ grows the period of
time for which they remain close to $W^{s}(y)$ grows arbitrarily.
Denote by $z_{n}$ the intersections of $W^{s}(y_{n})$ with
$W^{u}(f(y))$. $f(y)$ divides the unstable curve $W^{u}(f(y))$ in
a bounded arc $f(x)f(y)$ and an unbounded arc. We will prove that
$z_{n}$ belongs to the compact arc $f(x)f(y)$ for each $n$. If
some $z_{n}$ belongs to the unbounded arc, let us consider the
closed curve $J$ determined by the stable arc $xy_{n}$, the stable
arc $y_{n}z_{n}$, the unstable arc $z_{n}f(x)$ and the stable arc
$f(x)x$ (see figure \ref{fig231}).
\begin{figure}[htb]
\psfrag{p}{$p$} \psfrag{Wup}{$W^{u}(p)$} \psfrag{Wsp}{$W^{s}(p)$}
\psfrag{Wuy}{$W^{u}(y)$} \psfrag{Wsy}{$W^{s}(y)$}
\psfrag{Wsyn}{$W^{s}(y_{n})$} \psfrag{Wufy}{$W^{u}(f(y))$}
\psfrag{f(y)}{$f(y)$} \psfrag{y}{$y$} \psfrag{x}{$x$}
\psfrag{yn}{$y_{n}$}\psfrag{zn}{$z_{n}$} \psfrag{k}{$k$}
\psfrag{fx}{$f(x)$}
\begin{center}
\includegraphics[scale=0.25]{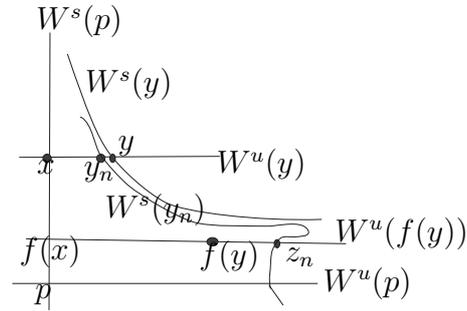}
\caption{\label{fig231}Invariant stable IV}
\end{center}
\end{figure}

Similar arguments as those used in the previous case allow us to
state that $W^{s}(f(y))$ can not intersect $J$ twice, which is a
contradiction. This proves that $z_{n}$ is in the bounded arc
$f(x)f(y)$, for every $n$. Let us consider, as figure \ref{fig74}
shows, the segment $y_{n}z_{n}$ and let $w_{n}$ be the point of
$W^{s}(y_{n})$ which is farthest from the fixed point.
\begin{figure}[htb]
\psfrag{p}{$p$} \psfrag{Wuy}{$W^{u}(y)$} \psfrag{Wsy}{$W^{s}(y)$}
\psfrag{Wsyn}{$W^{s}(y_{n})$} \psfrag{Wufy}{$W^{u}(f(y))$}
\psfrag{f(y)}{$f(y)$} \psfrag{y}{$y$}
\psfrag{yn}{$y_{n}$}\psfrag{zn}{$z_{n}$}
\psfrag{wn}{$w_{n}$}\psfrag{rn}{$r_{n}$}
\psfrag{Wuwn}{$W^{u}(w_{n})$}\psfrag{qn}{$q_{n}$}
\begin{center}
\includegraphics[scale=0.25]{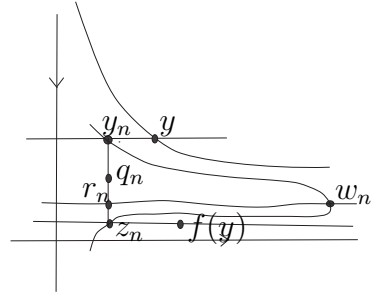}
\caption{\label{fig74}Invariant stable V}
\end{center}
\end{figure}
$W^{u}(w_{n})$ must intersect segment $y_{n}z_{n}$. Otherwise, it
would cut $W^{s}(y_{n})$ more than once. Let $r_{n}$ be that
intersection point. We want to apply our condition {\bf HP}.
Observe that points $r_{n},y_{n}$ would be in a compact set and
$w_{n}$ tends to infinity when $y_{n}$ tends to $y$.
$V(w_{n},r_{n})>0$ because they are in the same unstable curve,
and $V(w_{n},y_{n})<0$ because they are on the same stable curve.
So, there exists a point $q_{n}$ that belongs to segment
$y_{n}r_{n}$ such that $V(w_{n},q_{n})=0$. Then $$ \lim _{ n
\rightarrow \infty }
\frac{|V(w_{n},y_{n})-V(w_{n},q_{n})|}{W(w_{n},y_{n})} =0.$$
Therefore, we can choose $w_{n}$ such that
$$W(w_{n},y_{n})+V(w_{n},y_{n})>0,$$ which implies that
$$V(f(w_{n}),f(y_{n}))>0.$$ This contradicts the fact that points
$f(w_{n}),f(y_{n})$ are in the same stable curve.
\end{itemize}
Thus the existence of our invariant stable curve is proved. From
now on, we will denote it by $J$. Next, we will prove, using again
condition {\bf HP}, that the existence of this invariant stable
curve is not possible, so we will end the proof. Let us take any
point $x$ in the stable curve of the fixed point. We state that
the unstable curve $W^{u}(x)$ through $x$ intersects $J$.
Otherwise, there would exist a point $x_{0}$ in the stable curve
of the fixed point such that its unstable curve, $W^{u}(x_{0})$,
is the first one that does not intersect $J$. But then
$W^{u}(f^{-1}(x_{0}))$ would intersect $J$. This is a
contradiction since it is one of the previous cases. Then, we have
that every point $x$ of the stable curve of the fixed point has
the property that its unstable curve intersects $J$. Moreover, as
point $x$ comes closer to the fixed point, the intersection, $z$,
gets closer to infinity, since $J$ separates the plane. We want to
apply our hypothesis {\bf HP}. Let us fix a point $y$ at the
invariant stable curve $J$ and a point $s$ in the unstable curve
of the fixed point (see figure \ref{fig78}).
\begin{figure}[htb]
\psfrag{p}{$p$} \psfrag{Wup}{$W^{u}(p)$} \psfrag{Wsp}{$W^{s}(p)$}
\psfrag{J}{$J$} \psfrag{y}{$y$} \psfrag{x}{$x$} \psfrag{s}{$s$}
\psfrag{t}{$t$} \psfrag{zt}{$z_{t}$} \psfrag{qt}{$q_{t}$}
\begin{center}
\includegraphics[scale=0.25]{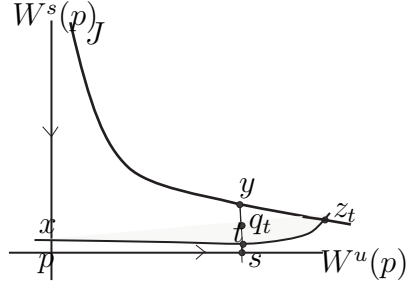}
\caption{\label{fig78}Final argument}
\end{center}
\end{figure}
Then, let us consider $x$ close enough to the fixed point, and let
$t$ be the intersection of $W^{u}(x)$ with segment $ys$. Let
$z_{t}$ be the intersection of $W^{u}(x)=W^{u}(t)$ with $J$.
Reasoning in a similar way to other parts of this proof, we have
that $V(z_{t},t)>0$ because they are in the same unstable curve,
and $V(z_{t},y)<0$ because they are in the same stable curve. So,
there exists a point $q_{t}$ in segment $yt$ such that
$V(z_{t},q_{t})=0$. If $x$ gets closer to $p$, $z_{t}$ tends to
infinity, and then we are able to choose $x$ such that
$$\frac{|V(z_{t},y)-V(z_{t},q_{t})|}{W(z_{t},y)} <1,$$ which implies
that $$W(z_{t},y)+ V(z_{t},y)>0, $$ and then
$$V(f(z_{t}),f(y))>0.$$ This yields a contradiction since points
$f(z_{t}),f(y)$ are in the same stable curve.\ep
\begin{lemma}
\label{LG1} Let $f$ be a homeomorphism of the plane that verifies
all the conditions described at the beginning of this section.
Then the stable (unstable) curve of every point intersects the
unstable (stable) curve of the fixed point.
\end{lemma}
\bp Let us consider set $A$ consisting of the points whose stable
(unstable) curve intersects the (unstable) stable curve of the
fixed point. It is clear that $A$ is open. Let us prove that it is
also closed. Let $(q_{n})$ be a sequence of $A$, convergent to
some point $q$ (see figure \ref{fig11}). Let $V(q)$ be a
neighborhood of $q$ with local product structure.
\begin{figure}[htb]
\psfrag{p}{$p$} \psfrag{Wuq}{$W^{u}(q)$} \psfrag{Wsq}{$W^{s}(q)$}
\psfrag{Wuqn0}{$W^{u}(q_{n_{0}})$}
\psfrag{Wsqn0}{$W^{s}(q_{n_{0}})$} \psfrag{V(q)}{$V(q)$}
\psfrag{q}{$q$} \psfrag{qn0}{$q_{n_{0}}$}
\psfrag{alphan0}{$\alpha_{n_{0}}$}
\begin{center}
\includegraphics[scale=0.25]{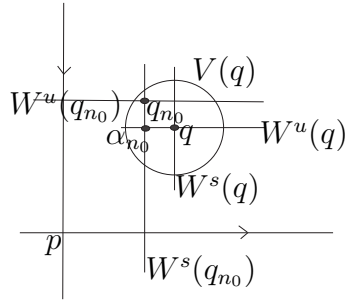}
\caption{\label{fig11}Coordinates}
\end{center}
\end{figure}
Let us consider $q_{n_{0}}\in V(q)$. So, we have that
$W^{s}(q_{n_{0}})\cap W^{u}(q)=\alpha_{n_{0}} $ as a consequence
of the local product structure and $W^{s}(q_{n_{0}})\cap
W^{u}(p)\neq \emptyset $ since $q_{n_{0}}\in A$. But then
$\alpha_{n_{0}} $ is a point in $W^{s}(q_{n_{0}})$ that cuts the
unstable curve of the fixed point, and then, applying lemma
\ref{TG} we have that $W^{u}(q)=W^{u}(\alpha_{n_{0}} )$ must cut
the stable curve of the fixed point. A similar argument lets us
prove that the stable curve of $q$ must cut the unstable curve of
the fixed point. Therefore $q$ belongs to set $A$ and consequently
$A$ is closed. Then $A$ is the whole plane.\ep
\begin{obs}
At this point we can not ensure that every stable (unstable) curve
cuts every unstable (stable) curve. Theorem \ref{TP}, one of the
main results in this work, shows that under our conditions, which
means admitting the existence of a Lyapunov metric function with
the required hypothesis, every stable (unstable) curve cuts every
unstable (stable) curve.
\end{obs}
\subsection{Main result.}
\begin{figure}[htb]
\psfrag{p1}{$p_{1}$} \psfrag{qx}{$q_{x}$} \psfrag{p2}{$p_{2}$}
\psfrag{p}{$p$} \psfrag{x}{$x$} \psfrag{Wsp1}{$W^{s}(p_{1})$}
\psfrag{Wup2}{$W^{u}(p_{2})$} \psfrag{zx}{$z_{x}$}
\begin{center}
\includegraphics[scale=0.3]{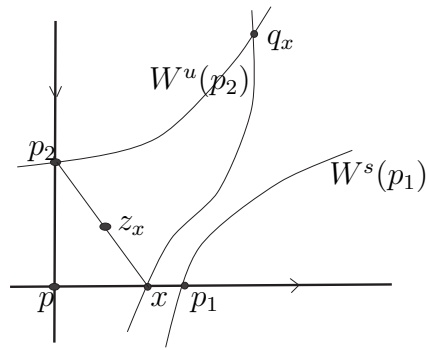}
\caption{\label{fig56}Argument.}
\end{center}
\end{figure}

\begin{figure}[htb]
\psfrag{H(q)}{$H(q)$} \psfrag{q}{$q$} \psfrag{H(x)}{$H(x)$}
\psfrag{x}{$x$} \psfrag{H(y)}{$H(y)$} \psfrag{y}{$y$}
\psfrag{H}{$H$}
\begin{center}
\includegraphics[scale=0.2]{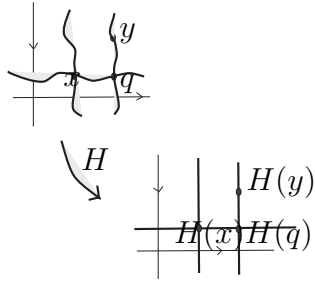}
\caption{\label{fig80}Argument.}
\end{center}
\end{figure}

In this section we will prove one of the main results of this work.
We obtain a characterization theorem (Theorem \ref{TP}) of expansive
homeomorphisms which verify the set of conditions exposed at the
beginning.
\begin{proposition}
Let $f:\R^{2} \rightarrow \R^{2}$ be a homeomorphism such that:
\begin{itemize}
\item $f$ has a fixed point; \item the quadrants determined by the
stable and unstable curves of the fixed point are $f$-invariant;
\item $f$ admits a Lyapunov metric function $U:\R^{2} \times
\R^{2} \rightarrow \R $. $U$ induces on the plane the same
topology than the usual distance; \item $f$ does not have
singularities; \item for each point $x\in \R^{2} $ and any $k>0$
there exist points $y$ and $z$ in the border of $B_{k}(x)$ such
that $V(x,y)=U(f(x),f(y))-U(x,y)>0$ and
$V(x,z)=U(f(x),f(z))-U(x,z)<0$, where $B_{k}(x)=\{ y\in \R^{2} \ /
U(x,y)\leq k\} $.
\end{itemize}
{\bf Then, if $U$ admits the condition (HP), then $f$ is
conjugated to a linear hyperbolic automorphism.}
\end{proposition}
\bp Let $f$ be a homeomorphism of the plane that admits a Lyapunov
metric function with the given hypothesis. If we proved that every
stable curve intersects every unstable curve, we would be able to
define a conjugation $H$ between a linear hyperbolic automorphism
$F$ and $f$, in the following way: first, we define $H$ sending
the stable (unstable) curve of the fixed point of $f$ to the
stable (unstable) curve of the fixed point of $F$ and such that
$F\circ H=H\circ f$. Since every $q$ on the plane is determined by
$\{ q\} =W^{s}(x)\cap W^{u}(y)$ with $x,y$ belonging respectively
to the unstable and stable curves of the fixed point, we define
$H(q)=W^{s}(H(x))\cap W^{u}(H(y))$. If we prove that every stable
curve intersects every unstable curve, we could say that $H$ is a
homeomorphism of $\R^{2} $ over $\R^{2} $. Let us prove that every
stable curve intersects every unstable curve. Let us suppose that,
as shown in figure \ref{fig56}, there exist $p_{1}\in W^{u}(p)$
and $p_{2}\in W^{s}(p)$ ($p$ is the fixed point) such that
$$W^{u}(p_{2})\cap W^{s}(p_{1})=\emptyset .$$ We are also under
the assumption that $p_{1}$ is the first point in $W^{u}(p)$ such
that its stable curve does not intersect $W^{u}(p_{2})$.
Thus $$W^{u}(p_{2})\cap W^{s}(x)=\{ q_{x}\} ,$$ were $x$ is a
point in the unstable arc $pp_{1}$ and $q_{x}$ near infinity as we
want, when $x$ tends to $p_{1}$. Let $V$ be the first difference
of the Lyapunov function. As points $x,p_{2}$ are in a compact
set, we are able to apply condition {\bf HP} concerning the
Lyapunov function in the following way: $V(q_{x},p_{2})>0$, since
both points are in the same unstable curve and $V(q_{x},x)<0$,
since both points are in the same stable curve. So, there exists a
point $z$ on segment $xp_{2}$ such that $V(z_{x},q_{x})=0$. Then
$$ \lim _{\| x\| \rightarrow p_{1} }
\frac{|V(q_{x},x)-V(q_{x},z_{x})|}{W(q_{x},x)} =0.$$ Then, we can
choose $q_{x}$ such that $$W(q_{x},x)+V(q_{x},x)>0,$$ which
implies that $$V(f(q_{x}),f(x))>0.$$ This contradicts the fact
that points $f(q_{x})$ and $f(x)$ are in the same stable
curve.\ep
\begin{theorem}
\label{TP} In the same conditions we had in the previous
proposition, $f$ is conjugated to a linear hyperbolic automorphism
if and only if it admits a Lyapunov metric function satisfying
condition {\bf (HP)}.
\end{theorem}
\bp It is a consequence of the last proposition and section
\ref{PRE}.\ep
\begin{corollary}
Under the same conditions of the previous theorem, $f$ admits a
Lyapunov function satisfying condition {\bf (HP)} if and only if
it is conjugated to a lift of an expansive homeomorphism on
$T^{2}$.
\end{corollary}
The stable curve $W^{s}$ and the unstable curve $W^{u}$ which we
built in this article, have the property that given two points
$x,y$ of $W^{s}$ ($W^{u}$) there exists $k>0$ such that
$U(f^{n}(x),f^{n}(y))<k$ for $n>0$ ($U(f^{n}(x),f^{n}(y))<k$ for
$n<0$). Denote by $W^{s}_{d}$ ($W^{u}_{d}$) the stable (unstable)
curves in the usual metric $d$ sense, this means that verifies
that given two points $x,y$ of $W^{s}_{d}$ ($W^{u}_{d}$) there
exists $k>0$ such that $d(f^{n}(x),f^{n}(y))<k$ for $n>0$
($d(f^{n}(x),f^{n}(y))<k$ for $n<0$). Observe that the
conjugations that appear in Proposition \ref{TP} send stable
(unstable) curves $W^{s}_{d}$ ($W^{u}_{d}$) of the linear
automorphism into stable (unstable) curves $W^{s}$ ($W^{u}$) of
$f$. The following proposition is a necessary and sufficient
condition to preserve the stable curves and unstable curves in the
sense of the usual metric.
\begin{proposition}
A necessary and sufficient condition for the conjugation of the
theorem \ref{TP} to preserve stable and unstable curves (in the
sense of the usual distance of the plane) is that the
homeomorphism $f$ admits a Lyapunov function $U$ that verifies the
following property: given any $k>0$ there exists $k'>0$ such that
$U(x,y)<k$ implies $d(x,y)<k'$, for each $x,y$ in $\R^{2} $.
\end{proposition}
\bp Let us consider a function $f$ that admits a Lyapunov function
$U$ with the property of the statement. Because of Theorem
\ref{TP}, $f$ is conjugated to a linear hyperbolic automorphism
and the conjugacy $H$ maps linear stable (unstable) curves into
stable (unstable) curves of $f$ in the sense of function $U$. But,
precisely this stable (unstable) curve satisfies
$U(f^{n}(x),f^{n}(y))<k$ for some $k>0$ and every $n>0$ ($n<0$).
Then, using the property, we have that there exists $k'>0$ such
that $d(f^{n}(x),f^{n}(y))<k'$, which proves that $H$ preserves
stable and unstable curves in the sense of the usual metric $d$.
Let us suppose now that $H$ preserves stable and unstable curves
in the sense of the usual metric. Let us consider two points $x,y$
such that $U(x,y)<k$, where $U(x,y)=D(H(x),H(y))$
($D=D_{s}+D_{u}$). Let $q$ be the intersection of the unstable
curve of $x$ with the stable curve of $y$ (see figure
\ref{fig80}).
As $U(x,y)<k$ implies $D(H(x),H(y))<k$, then
$$D_{s}(H(x),H(q))=D(H(x),H(q))<k$$ and
$$D_{u}(H(y),H(q))=D(H(y),H(q))<k.$$ Since $H$ preserves stable and
unstable curves, $$d(x,q)<k_{1}$$ and $$d(y,q)<k_{1},$$ with a
uniform $k_{1}$. This implies that there exists $k_{2}$ uniform
such that $d(x,y)<k_{2}$ which concludes the proof.\ep
\begin{corollary}
A homeomorphism of the plane with the conditions given in this
section is the time $1$ of a flow.
\end{corollary}
\bp We have that $$f=H^{-1}\varphi_{1} H,$$ where $\varphi_{1}$ is
a linear automorphism. Then we can consider the flow
$$\psi_{t} =H^{-1}\varphi_{t} H$$ and this implies that $f$ is $\psi_{1}$.\ep

\section{Another context and some examples.}
\label{OBS} In this section we will show some generalizations
about the hypothesis we asked for our homeomorphisms. The main
difference with what we have exposed until now, is that we will
work without condition {\bf HP}. Instead of this we will ask for
uniform local conditions concerning the first and the second
difference ($V$ and $W$) of the Lyapunov function $U$. In this new
context, we will get a new characterization result in Theorem
\ref{LOC}, which shows two possible behaviors of our
homeomorphisms. At the end of this section we will show some
examples.

\subsection{Local properties.} \label{SPL}
We will denote the following hypothesis for a homeomorphism $f$ by
condition {\bf HL}:
\begin{enumerate}
\item $f$ has a fixed point; \item the quadrants determined by the
stable and unstable curves of the fixed point are $f$-invariant;
\item $f$ admits a metric Lyapunov function $U:\R^{2} \times
\R^{2} \rightarrow \R $. $U$ induces in the plane the same
topology as the usual distance; \item for each point $x\in \R^{2}
$ and any $k>0$ there exist points $y$ and $z$ in the border of
$B_{k}(x)$ such that $V(x,y)=U(f(x),f(y))-U(x,y)>0$ and
$V(x,z)=U(f(x),f(z))-U(x,z)<0$, where $B_{k}(x)=\{ y\in \R^{2} \ /
U(x,y)\leq k\} $; \item $f$ does not admit singularities; \item
\label{ARC} given any $\epsilon >0$ and two points $x,y$ of
$\R^{2}$ such that $U(x,y)<\epsilon $, there exists an arc $a$
that joins $x$ with $y$ such that $U(z,t)<\epsilon $ for each pair
of points $z,t$ that belong to arc $a$; \item \label{k} given any
$k>0$ we denote by $W^{u}_{k}(x)$ the $k$-unstable arc of $x$ i.e.
the connected component of the set $\{y\in \R^{2}:
U(f^{n}(x),f^{n}(y))\leq k, \ \forall n\leq 0\}$ that contains
$x$. Let $x,y\in \R^{2}$ be such that $U(f^{n}(x),f^{n}(y)$ tends
to zero when $n$ tends to infinity. There exists $k(x,y)>0$ such
that $U(z_{n},t_{n})$ tends to zero, where $z_{n},t_{n}$ are
endpoints of $W^{u}_{k}(f^{n}(x))$ and $W^{u}_{k}(f^{n}(y))$
respectively (see figure \ref{fig400});

\begin{figure}[htb]
\psfrag{fnx}{$f^{n}(x)$} \psfrag{fny}{$f^{n}(y)$}
\psfrag{zn}{$z_{n}$} \psfrag{tn}{$t_{n}$}
\psfrag{Wukfnx}{$W^{u}_{k}(f^{n}(x))$}
\psfrag{Wukfny}{$W^{u}_{k}(f^{n}(y))$}
\begin{center}
\includegraphics[scale=0.2]{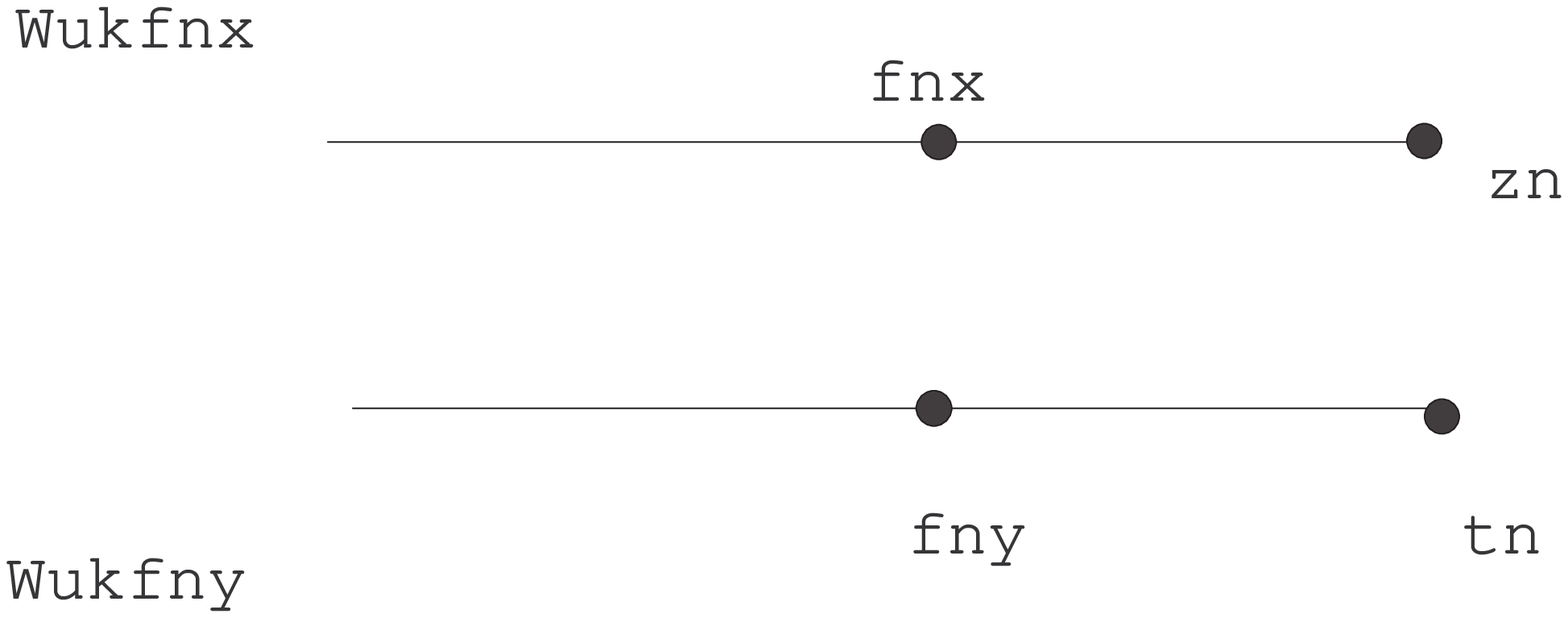}
\caption{\label{fig400}Hypothesis.}
\end{center}
\end{figure}
\item \label{V} The first difference $V=\Delta U$ verifies the
following property: given any $\epsilon >0$ there exists $\delta
>0$ such that if $U(z,y)<\delta $ then $|V(x,z)-V(x,y)|<\epsilon
$, for all $x\in \R^{2}$; \item \label{W} The second difference
$W=\Delta^{2} U$ verifies the following property: given any
$\delta >0$ there exists $a(\delta )>0$ such that $W(x,y)>a(\delta
)>0$ for all $x,y$ in the plane with $U(x,y)>\delta $.
\end{enumerate}

\begin{lemma}
\label{LPL} Let $f:\R^{2} \rightarrow \R^{2}$ be a homeomorphism
of the plane that verifies condition {\bf HL}. If two points $x,y$
verify that $U(f^{n}(x),f^{n}(y))\leq U(x,y)$ for all $n>0$, then
they belong to the same stable curve. Moreover $U(f^{n}(x),
f^{n}(y))$ tends to zero, when $n$ tends to infinity. A similar
statement can be given for $n\leq 0$.
\end{lemma}

\bp Let us consider two points $x,y$ such that
$U(f^{n}(x),f^{n}(y))\leq U(x,y)$ for all $n>0$. So,
$V(f^{n}(x),f^{n}(y))<0$ for each $n>0$, because if for some $n>0$
$V(f^{n}(x),f^{n}(y))>0$, we would have that
$V(f^{m}(x),f^{m}(y))>V(f^{n}(x),f^{n}(y))$ if $m>n$ ($W>0$) and
then $$U(f^{m}(x),f^{m}(y))=U(f^{n}(x),f^{n}(y))+\sum_{j=n}^{m}
V(f^{j}(x),f^{j}(y))>$$
$$U(f^{n}(x),f^{n}(y))+(m-n)V(f^{n}(x),f^{n}(y))$$ which implies that
$U(f^{m}(x),f^{m}(y))$ is unbounded, which contradicts the
assumption. Then, $V(f^{n}(x),f^{n}(y))<0$ for $n>0$. Since $\sum
V(f^{n}(x),f^{n}(y))$ is bounded, $V(f^{n}(x),f^{n}(y))$ tends to
$0$ when $n$ tends to infinity. So, $U(f^{n}(x),f^{n}(y))$ tends
to $0$. Otherwise, there would exist $\delta >0$ and $n>0$ large
enough such that $U(f^{n}(x), f^{n}(y))>\delta$, and then
$W(f^{n}(x),f^{n}(y))>a(\delta)>0$ (property \ref{W}). As
$V(f^{n}(x),f^{n}(y))$ tends to $0$ when $n$ tends to infinity and
$$W(f^{n}(x),f^{n}(y))= V(f^{n+1}(x),f^{n+1}(y))-V(f^{n}(x),f^{n}(y)),$$
we would find some $n>0$ such that $V(f^{n+1}(x),f^{n+1}(y))>0$
which is not possible. Then, $U(f^{n}(x), f^{n}(y))$ tends to $0$
when $n$ tends to infinity. Let $u_{n}=f^{n}(x), v_{n}=f^{n}(y)$
and let $p_{n}$ and $q_{n}$ be the $k$- unstable arc through
$u_{n}$ and $v_{n}$ respectively, with $k$ given by hypothesis
\ref{k} of condition {\bf HL}. Now, we will prove that points
$u_{n}$ and $v_{n}$ would be in the same global stable curve. We
will divide the reasoning in two cases: \begin{itemize} \item {\it
The stable curve through $v_{n}$ intersects $p_{n}$ or the stable
curve through $u_{n}$ intersects $q_{n}$.} Let us suppose (as
shown in figure \ref{fig61}) $w_{n}=W^{s}(v_{n})\cap p_{n}$,
$w_{n}\neq u_{n}$.
\begin{figure}[htb]
\psfrag{pn}{$p_{n}$} \psfrag{qn}{$q_{n}$} \psfrag{un}{$u_{n}$}
\psfrag{wn}{$w_{n}$} \psfrag{vn}{$v_{n}$}
\psfrag{Wsvn}{$W^{s}(v_{n})$}
\begin{center}
\includegraphics[scale=0.3]{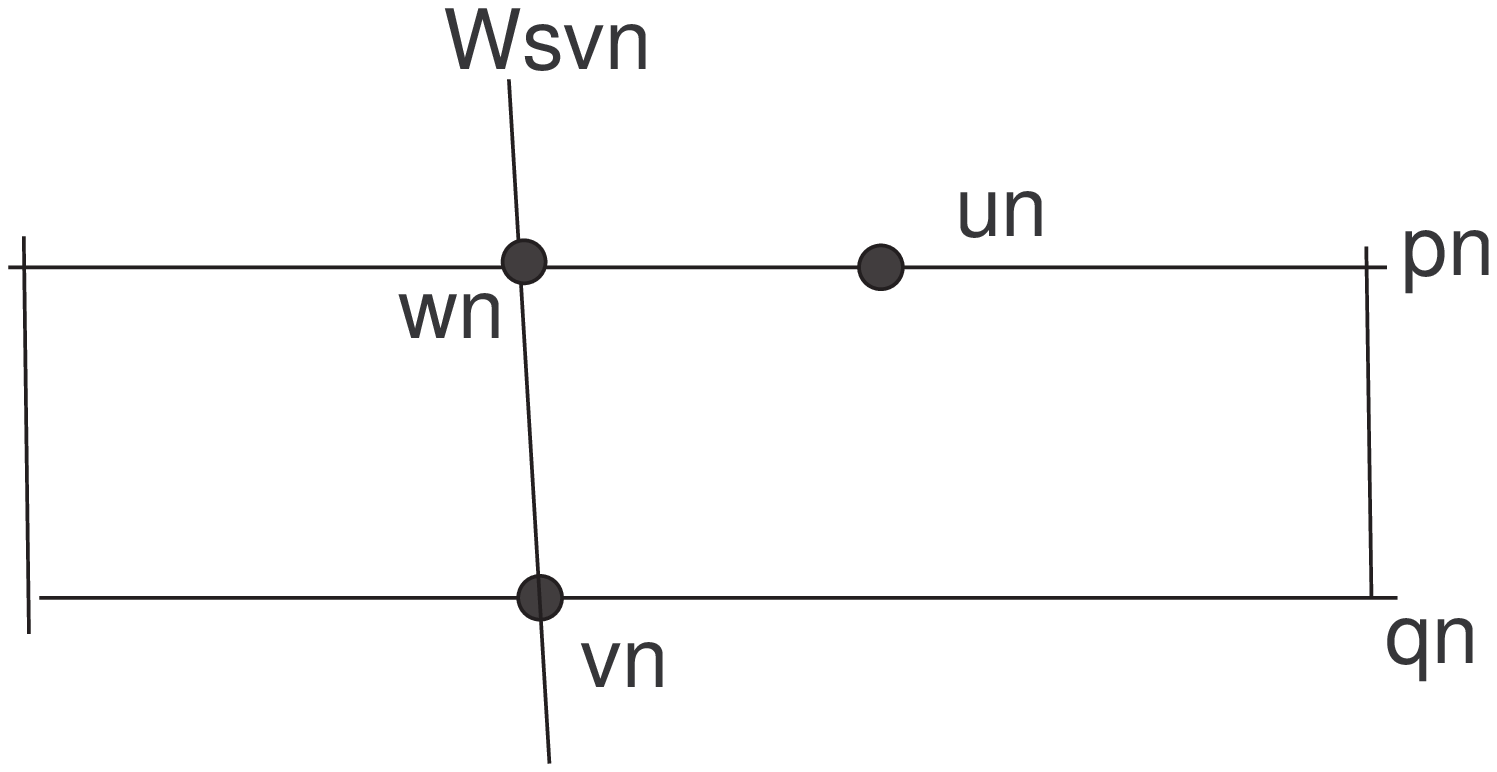}
\caption{\label{fig61}Case I.}
\end{center}
\end{figure}
It has already been proved that $U(u_{n},v_{n})$ and
$V(u_{n},v_{n})$ tends to zero when $n$ tends to infinity. Then,
using hypothesis \ref{V} for $V$, we have that
$|V(u_{n},w_{n})-V(v_{n},w_{n})|$ must tend to zero when $n$ tends
to infinity. But $V(u_{n},w_{n})$ is positive and grows with $n$
(because $u_{n}$ and $v_{n}$ are in the same unstable curve),
which yields a contradiction. Then $w_{n}=u_{n}$. \item {\it The
stable curve through $v_{n}$ does not intersect $p_{n}$ and the
stable curve through $u_{n}$ does not intersect $q_{n}$.} As
$U(u_{n},v_{n})$ tends to zero when $n$ tends to infinity then,
using hypothesis \ref{k} of condition {\bf HL}, we have that
$U(\alpha_{n} ,\beta_{n} )$ tends to zero, where $\alpha_{n}
,\beta_{n} $ are endpoints of $p_{n}$ and $q_{n}$ (see figure
\ref{fig54}).
\begin{figure}[htb]
\psfrag{pn}{$p_{n}$} \psfrag{qn}{$q_{n}$} \psfrag{un}{$u_{n}$}
\psfrag{wn}{$w_{n}$} \psfrag{vn}{$v_{n}$}
\psfrag{alphan}{$\alpha_{n}$} \psfrag{betan}{$\beta_{n}$}
\psfrag{an}{$a_{n}$}
\begin{center}
\includegraphics[scale=0.3]{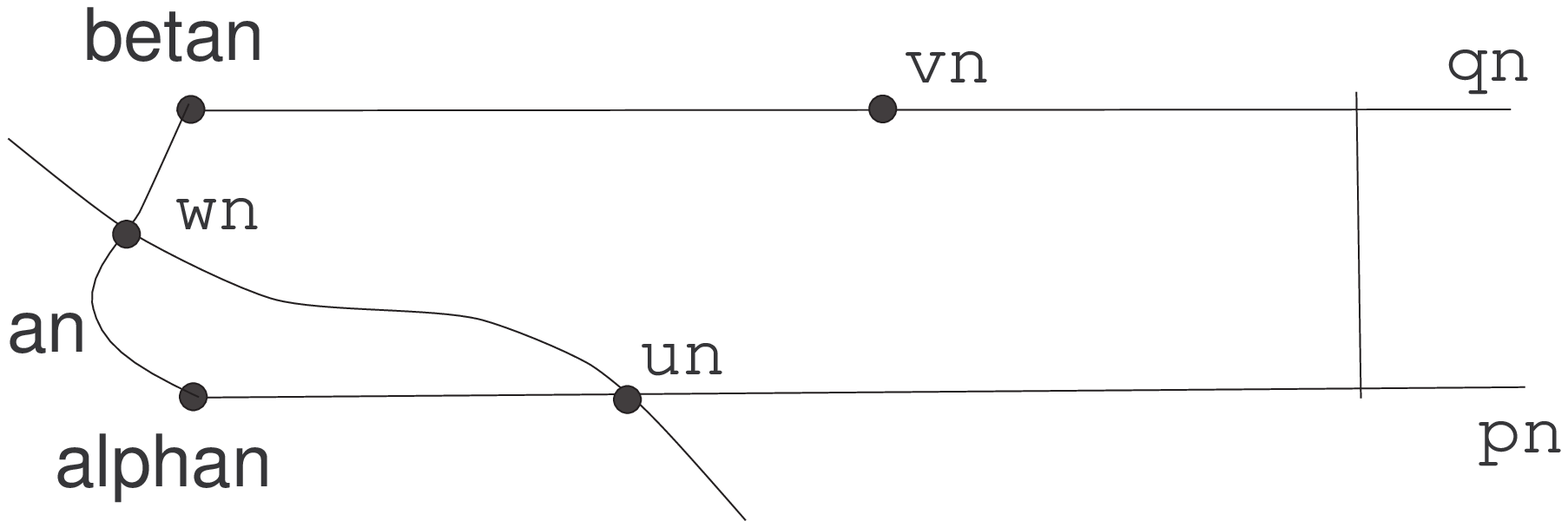}
\caption{\label{fig54}Caso II.}
\end{center}
\end{figure}
Using hypothesis \ref{ARC} of condition {\bf HL}, there exists an
arc $a_{n}$ that joins points $\alpha_{n} $ and $\beta_{n} $ such
that if $U(\alpha_{n} ,\beta_{n} )<\epsilon $ then $U(\alpha_{n}
,w_{n})<\epsilon $ for each $w_{n}$ on arc $a_{n}$. If the stable
curve through $u_{n}$ does not intersect $q_{n}$, it has to
intersect one of the arcs that join the endpoints of $p_{n}$ and
$q_{n}$ (recall that since the considered stable curve is global,
it separates the plane). Let us suppose, without loosing
generality, that this arc is $a_{n}$. Let $\{w_{n}\}
=W^{s}(u_{n})\cap a_{n}$. Since $$U(u_{n},\alpha_{n})\leq
U(\alpha_{n},w_{n})+U(u_{n} ,w_{n})$$ and
$$U(u_{n},\alpha_{n})\geq \delta (k)>0,$$ we have that $U(u_{n}
,w_{n})$ is bounded away from zero for each $n>0$. Therefore,
applying hypothesis \ref{W} we have that
$W(u_{n},w_{n})>k_{1}(k)>0$. $V(u_{n},w_{n})<0$ since these points
are in the same stable curve and $V(u_{n},\alpha_{n})>0$ since
these points are in the same unstable curve. As $U(\alpha_{n}
,w_{n})$ tends to zero, then we can apply hypothesis \ref{V} and
so we can state that $V(u_{n},w_{n})$ is arbitrarily close to zero
for some $n$. Then, $$V(f(u_{n}),f(w_{n}))=W(u_{n},w_{n})+
V(u_{n},w_{n})>0,$$ which contradicts the fact that $f(u_{n})$ and
$f(w_{n})$ are in the same stable curve.\end{itemize}
\ep

\begin{lemma}
\label{TL} Let $f:\R^{2} \rightarrow \R^{2}$ be a homeomorphism of
the plane that verifies condition {\bf HL}. Let us suppose that
there exists a point $x$ of the stable curve $W^{s}(p)$ of the
fixed point, such that there exists a point $z$ belonging to
$W^{u}(x)$ that verifies that $W^{s}(z)\cap W^{u}(p)=\emptyset $.
Then, there exists a point $y\in W^{u}(x)$ such that $W^{s}(y)$ is
invariant under $f$.
\end{lemma}
\bp The first part of this proof uses arguments which are similar
to those that we used in Lemma \ref{TG}. Let us remember the
beginning of the proof of this lemma: let $y$ be the first point
of $W^{u}(x)$ with the property that its stable curve does not
intersect the unstable curve of the fixed point (see figure
\ref{fig13}).
\begin{figure}[htb]
\psfrag{p}{$p$} \psfrag{Wsy}{$W^{s}(y)$} \psfrag{Wux}{$W^{u}(x)$}
\psfrag{Wsfy}{$W^{s}(f(y))$} \psfrag{f(y)}{$f(y)$} \psfrag{y}{$y$}
\psfrag{x}{$x$} \psfrag{B}{$B$}\psfrag{f-1B}{$f^{-1}(B)$}
\psfrag{(I)}{$(I)$} \psfrag{(II)}{$(II)$} \psfrag{Wsp}{$W^{s}(p)$}
\begin{center}
\includegraphics[scale=0.25]{fig13.eps}
\caption{\label{fig13}}
\end{center}
\end{figure}
We show that $W^{s}(y)$ is $f$ invariant. Let us reason in one
quadrant determined by the stable and the unstable curves of the
fixed point. Since $W^{s}(y)$ separates the first quadrant, denote
by zone $(I)$ the region that does not include the fixed point $p$
in its border, and zone $(II)$ the other one. We will divide the
proof in three steps:
\begin{itemize}
\item {\bf $f(y)$ belongs to the zone $(I)$} See Lemma \ref{TG} of
section \ref{SPG}. \item {\bf $f(y)$ belongs to the zone $(II)$
and $W^{u}(f(y))\cap W^{s}(y)= \{ q\} $.} See Lemma \ref{TG} of
section \ref{SPG}. \item {\bf $f(y)$ belongs to the zone $(II)$
and $W^{u}(f(y))\cap W^{s}(y)= \emptyset $.} Let us consider
$(y_{n})$, a sequence in $W^{u}(y)$ converging to $y$ such that
the stable curves of the points $y_{n}$ cut the unstable curve of
the fixed point. The behavior of the stable curves of $y_{n}$ must
be the one shown in figure \ref{fig12}: as $n$ grows the period of
time for which they remain close to $W^{s}(y)$ grows arbitrarily.
\begin{figure}[htb]
\psfrag{p}{$p$} \psfrag{Wup}{$W^{u}(p)$} \psfrag{Wsp}{$W^{s}(p)$}
\psfrag{Wuy}{$W^{u}(y)$} \psfrag{Wsy}{$W^{s}(y)$}
\psfrag{Wsyn}{$W^{s}(y_{n})$} \psfrag{Wufy}{$W^{u}(f(y))$}
\psfrag{f(y)}{$f(y)$} \psfrag{y}{$y$} \psfrag{x}{$x$}
\psfrag{yn}{$y_{n}$}\psfrag{zn}{$z_{n}$}
\begin{center}
\includegraphics[scale=0.25]{fig12.eps}
\caption{\label{fig12}}
\end{center}
\end{figure}
Denote by $z_{n}$ the intersections of $W^{s}(y_{n})$ with
$W^{u}(f(y))$. The first thing we will state is that sequence
$(z_{n})$ accumulates in a point $h$ of $W^{u}(f(y))$. Indeed,
$f(y)$ divides the unstable arc $W^{u}(f(y))$ in an bounded arc
$f(x)f(y)$ and in one unbounded arc. Using identical arguments to
those used in Lemma \ref{TG}, we can prove that $z_{n}$ belongs to
the bounded arc for each $n$ and this way we prove the existence
of an accumulation point $h$. Since points $h$ and $y$ are not on
the same stable curve, and considering the conclusion of Lemma
\ref{LPL}, we can state that $U(f^{n}(h),f^{n}(y))$ reaches
arbitrarily large values for certain $n>0$. Now, since
$y_{n}\rightarrow y$ and $z_{n}\rightarrow h$ we can state that
there exist $n,p\in \N $ such that
$U(f^{n}(z_{p}),f^{n}(y_{p}))>U(y,h)$. Denote by $M=U(y,h)$ and
choose $p\in \N $ in such a way that for some $n_{0}\in \N $,
$U(f^{n_{0}}(z_{p}),f^{n_{0}}(y_{p}))>>M$. Because of the
continuity of $U$ we have that for $p$ sufficiently large
$U(z_{p},y_{p})$ would be close to $M=U(y,h)$. This implies that
at some moment $U(f^{n}(z_{p}),f^{n}(y_{p}))$ grew, which means
that $V(f^{n}(z_{p}),f^{n}(y_{p}))>0$ for some $n\in \N $. Since
$W>0$, $U(f^{n}(z_{p}),f^{n}(y_{p}))$ will grow to infinity, which
contradicts the fact that points $z_{p}$ and $y_{p}$ are in the
same stable curve.\end{itemize}
\ep

The following condition will ensure us that we only have one
invariant stable curve: the stable curve of the fixed point.\\
{\bf Additional hypothesis.(HA)} {\it $f$ satisfies
$$\lim_{n\rightarrow \pm \infty } U(f^{n}(x),f^{n+1}(x))=\infty ,$$
for every $x\in \R^{2}$ that does not belong to the stable or
unstable curve of the fixed point.}\\ We will omit the proof of
the following lemma since it is analogous to Lemma \ref{LG1} of
section \ref{SPG}.

\begin{lemma}
Let $f$ be a homeomorphism of the plane that verifies condition
{\bf HL} and hypothesis {\bf HA}. Then, the stable (unstable)
curve of every point intersects the unstable (stable) curve of the
fixed point.
\end{lemma}

The following theorem will characterize homeomorphisms of the
plane that admit a Lyapunov function under condition {\bf HL} and
hypothesis {\bf HA}.

\begin{theorem}
\label{LOC} Let $f$ be a homeomorphism of the plane. Then, $f$
admits a Lyapunov function $U$ that verifies condition {\bf HL}
and hypothesis {\bf HA} if and only if $f$ restricted to each of
the quadrants determined by the stable and unstable curves of the
fixed point is either conjugated to a linear hyperbolic
automorphism or conjugated to the restriction of a linear
hyperbolic automorphism to certain invariant region. This
conjugations preserve stable and unstable curves.
\end{theorem}

\begin{figure}[htb]
\psfrag{h(p)}{$H(p)$} \psfrag{h(Wup)}{$H(W^{u}(p))$}
\psfrag{h(Wsp)}{$H(W^{s}(p))$} \psfrag{s}{$s$}
\begin{center}
\includegraphics[scale=0.2]{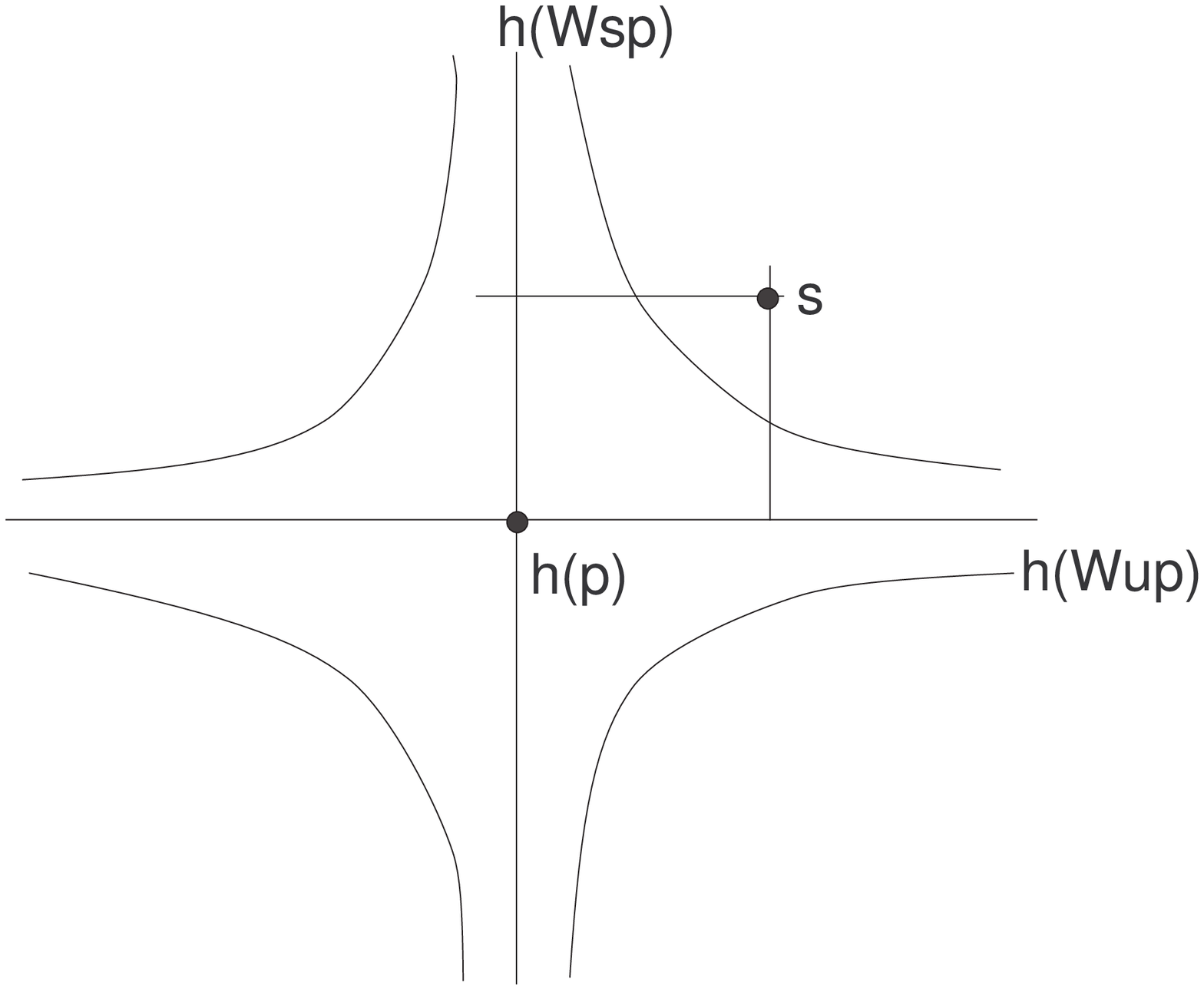}
\caption{\label{fig4}range}
\end{center}
\end{figure}

\bp Let us suppose that $f$ admits a Lyapunov function $U$ that
verifies condition {\bf HL} and hypothesis {\bf HA}. Let us try to
build the conjugation $H$ with the linear automorphism $F$. Let us
define $H$ sending the stable (unstable) curve of the fixed point
of $f$ on the stable (unstable) curve of the fixed point of $F$
and such that $F\circ H=H\circ f$. Then, given any point $q$ of
the selected quadrant we know, based on the previous results, that
$q=W^{s}(x)\cap W^{u}(y)$, with $x,y$ belonging to the unstable
and stable curve of the fixed point respectively. We define
$H(q)=W^{s}(H(x))\cap W^{u}(H(y))$. If the range of $H$ is the
whole quadrant, we will find a conjugation with the linear
automorphism on this quadrant. Otherwise, the range of $H$ is a
restriction of the linear automorphism to an invariant region
limited by the stable and unstable curves of the fixed point and a
decreasing curve, as shown in figure \ref{fig4}.

\begin{figure}[htb]
\psfrag{p}{$p$} \psfrag{a}{$a$}
\begin{center}
\includegraphics[scale=0.2]{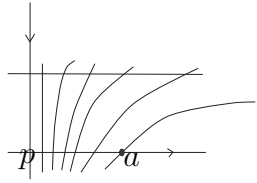}
\caption{\label{fig6}Behaviour}
\end{center}
\end{figure}

Notice that point $s$ does not have a preimage through $H$. This
case corresponds to one behavior of our homeomorphism as shown in
figure \ref{fig6}, where the stable curve of point $a$ does not
intersect the considered unstable curve. Later we will prove that
all these restrictions of the linear automorphism are conjugated
amongst themselves.

Conversely, let $f$ be a homeomorphism of the plane conjugated to
a linear automorphism $F$ or a restriction of this to an invariant
region. Define:

$$F(x,y)= \left(
\begin{array}{cc}
  \lambda & 0 \\
  0 & 1/\lambda \\
\end{array}
\right) \left( \begin{array}{c}
  x \\
  y \\
\end{array} \right) ,$$
with $\lambda >1$, let $D=D_{s}+D_{u}$ be the Lyapunov metric
function for $F$ (or its restriction) and define
$$L(p_{1},p_{2})=D(H^{-1}(p_{1}),H^{-1}(p_{2})).$$ Then, $L$ is a
Lyapunov metric function for $f$. The arguments are identical to
the ones used in Section \ref{PRE}. Now, we will verify that $L$
satisfies condition {\bf HL} and hypothesis {\bf HA}:
\begin{itemize}

\item {\it The first difference $\Delta L$ verifies the following
property: given any $\epsilon >0$ there exists $\delta >0$ such
that if $L(z,y)<\delta $ then $|\Delta L(x,z)-\Delta
L(x,y)|<\epsilon $.} It is easy to see that $$\Delta L(x,z)=\Delta
D(H^{-1}(x),H^{-1}(z))$$ where $$L(x,z)=D(H^{-1}(x),H^{-1}(z)).$$
In section \ref{PRE}, we showed that $$|\Delta
D(H^{-1}(x),H^{-1}(z))- \Delta D(H^{-1}(x),H(y))|\leq $$
$$(\lambda -1)D_{u}(H^{-1}(y),H^{-1}(z)) +$$ $$(1-1/\lambda
)D_{s}(H^{-1}(y),H^{-1}(z))\leq $$
$$KD(H^{-1}(y),H^{-1}(z))=KL(y,z),$$ which proves the property.

\item {\it The second difference $W=\Delta^{2} L$ verifies the
following property: given any $\delta >0$ there exists $a(\delta
)>0$ such that $W(x,y)>a(\delta )>0$ for each $x,y$ on the plane
with $L(x,y)>\delta $.} It is easy to see that $$\Delta^{2}
L(x,y)=\Delta^{2} D(H^{-1}(x),H^{-1}(y)).$$ In section \ref{PRE}
we showed that $$ \Delta^{2} D(x,y)= \Delta D(F(x),F(y))-\Delta
D(x,y)=$$
$$(\lambda -1)^{2}D_{u}(x,y) + (1-1/\lambda )^{2}D_{s}(x,y).$$ Then,
there exists $K>0$ such that
$$\Delta^{2} L(x,y)=\Delta^{2} D(H^{-1}(x),H^{-1}(y))\geq $$
$$KD(H^{-1}(x),H^{-1}(y))=KL(x,y).$$ Thus, this property is proved.

\item {\it For each point $x\in \R^{2} $ and any $k>0$ there exist
points $y$ and $z$ in the border of $B_{k}(x)$ such that $\Delta
L(x,y)=L(f(x),f(y))-L(x,y)>0$ and $\Delta
L(x,z)=L(f(x),f(z))-L(x,z)<0$.} The stable and unstable curves
separate the quadrant and that is why the property holds.

\item {\it Given any $\epsilon >0$ and two points $x,y$ of
$\R^{2}$ such that $L(x,y)<\epsilon $, there exists an arc $a$
that joins $x$ with $y$ such that $L(z,t)<\epsilon $ for each pair
of points $z,t$ that belongs to the arc $a$.} By definition,
$L(x,y)<\epsilon $ implies that $D(H^{-1}(x),H^{-1}(y))<\epsilon
$. Let us define segment $b$ as the one that joins points
$H^{-1}(x)$ and $H^{-1}(y)$. Then, arc $a=H(b)$ verifies the
property.

\item {\it Let $x,y\in \R^{2}$ be such that $L(f^{n}(x),f^{n}(y))$
tends to zero when $n$ tends to infinity. Then, there exists
$k(x,y)>0$ such that $L(z_{n},t_{n})$ tends to zero, where
$z_{n},t_{n}$ are endpoints of $W^{u}_{k}(f^{n}(x))$ and
$W^{u}_{k}(f^{n}(y))$, respectively.} Since $L(f^{n}(x),f^{n}(y)$
tends to zero, then $D(H^{-1}(f^{n}(x)),H^{-1}(f^{n}(y)))$ tends
to zero, which implies that $D(F^{n}(H^{-1}(x)),F^{n}(H^{-1}(y)))$
tends to zero. Then, points $H^{-1}(x)$ and $H^{-1}(y)$ are in the
same stable segment of the linear automorphism. Let $k>0$ be such
that the $k$-unstable segments of points $H^{-1}(x)$ and
$H^{-1}(y)$ belong to the invariant region considered. In the
future, the unstable segments with length $k$ of points
$F^{n}(H^{-1}(x))$ and $F^{n}(H^{-1}(y))$ will also be included in
the considered invariant region and its endpoints $H^{-1}(z_{n})$
and $H^{-1}(t_{n})$ will verify that
$D(H^{-1}(z_{n}),H^{-1}(t_{n}))$ (because they are in the same
stable segment of the linear automorphism) tends to zero. But then
$L(z_{n},t_{n})$ tends to zero as we wanted to prove.

\item {\it $f$ satisfies $$\lim_{n\rightarrow \pm \infty }
L(f^{n}(x),f^{n+1}(x))=\infty ,$$ for each $x\in \R^{2}$ that does
not belong to the stable or unstable curve of the fixed point.} We
know that
$$L(f^{n}(x),f^{n+1}(x))= D(H^{-1}(f^{n}(x)),H^{-1}(f^{n+1}(x)))=$$
$$D(F^{n}(H^{-1}(x)),F^{n+1}(H^{-1}(x))).$$ But
$D(F^{n}(H^{-1}(x)),F^{n+1}(H^{-1}(x)))$ tends to infinity when
$n$ tends to $\pm \infty $, since in the linear case there are no
invariant stable or unstable curve other than the stable or
unstable curves of the fixed point.\end{itemize}
\ep

\begin{proposition}
All the restrictions of the linear automorphism shown in Theorem
\ref{LOC} are conjugated amongst themselves.
\end{proposition}
\bp We will divide the proof in three cases:
\begin{itemize}

\item[Case 1] Let us suppose that we have two restrictions, $f$
and $g$, of the linear automorphism in the regions limited by
curves $J_{1}$ and $J_{2}$ as shown in figure \ref{fig17}. In
other words: there are no parts of these border curves that
consisting of segments which are parallel to the axis (stable and
unstable curves of the fixed point).

\begin{figure}[htb]
\psfrag{J1}{$J_{1}$} \psfrag{J2}{$J_{2}$} \psfrag{x}{$x$}
\psfrag{x1}{$x_{1}$} \psfrag{x2}{$x_{2}$} \psfrag{p1}{$p_{1}$}
\psfrag{p2}{$p_{2}$} \psfrag{H}{$H$} \psfrag{H(x1)}{$H(x_{1})$}
\psfrag{H(x2)}{$H(x_{2})$} \psfrag{H(x)}{$H(x)$}
\begin{center}
\includegraphics[scale=0.2]{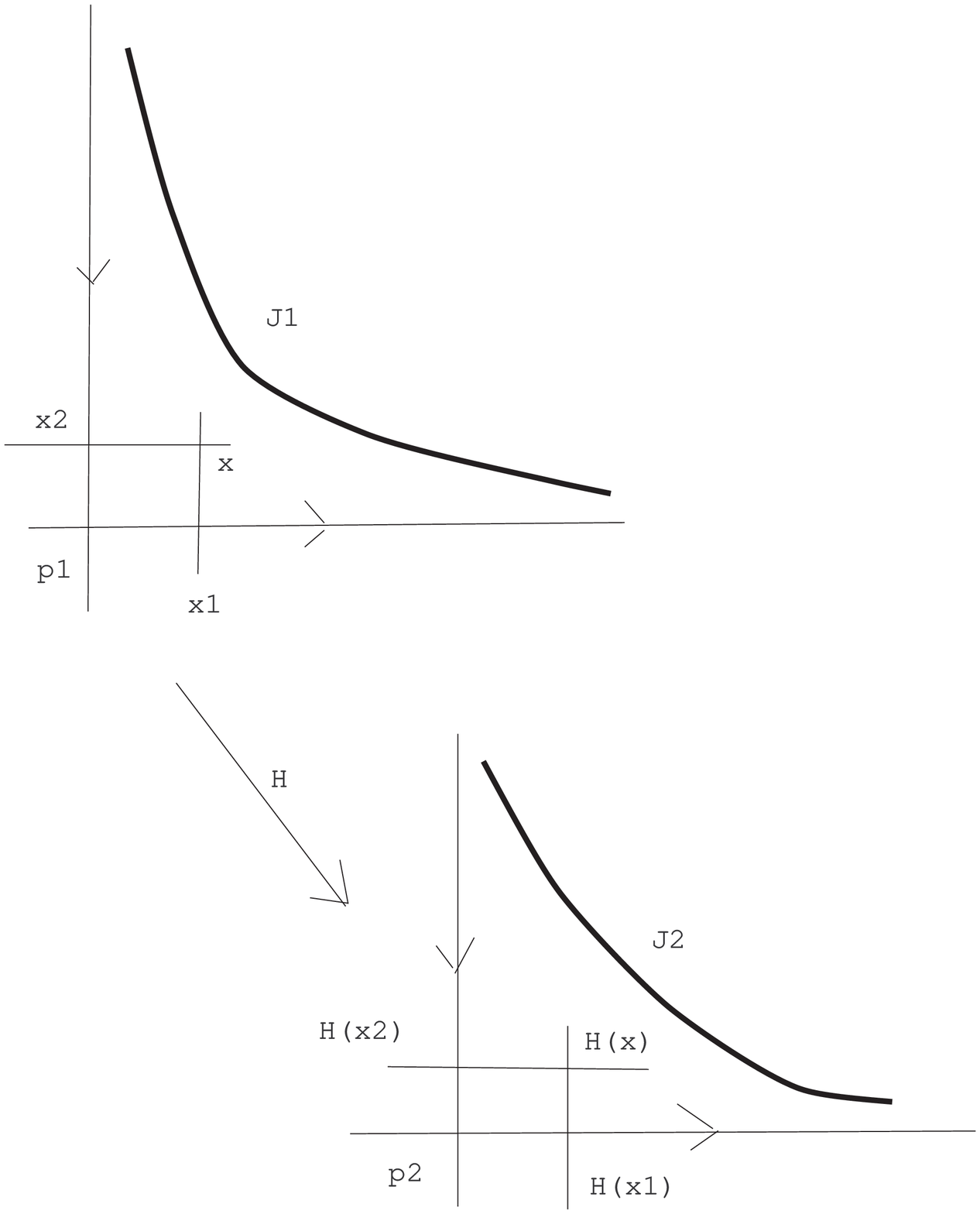}
\caption{\label{fig17}Case 1}
\end{center}
\end{figure}

Define $H$ mapping $J_{1}$ in $J_{2}$ such that $g\circ H=H\circ
f$. Now, we want to extend $H$ so that we can map the unstable
(stable) curve of the fixed point of $f$ in the unstable (stable)
curve of the fixed point of $g$. We will do it in the following
way: let $x\in W^{s}(p_{1})$. We define $$H(x)=W^{s}(p_{2})\cap
W^{u}(H(W^{u}(x)\cap J_{1})).$$ A similar definition for the
unstable curve. Then we are able to extend $H$ to the interior set
of the region of our interest (the one limited by $J_{1}$ and the
stable and unstable curves of the fixed point ($p_{1}$)) in the
following way: let $x$ be a point of this interior, so
$x=W^{s}(x_{1})\cap W^{u}(x_{2})$, where $x_{1}\in W^{u}(p_{1})$
and $x_{2}\in W^{s}(p_{1})$. We define
$$H(x)=W^{s}(H(x_{1}))\cap W^{u}(H(x_{2})).$$ This is the
conjugation we were searching for. \item[Case 2] Let us suppose
that we have some restriction of a lineal hyperbolic automorphism
to a region with border $J$, as the one shown in the figure
\ref{fig19}, that admits segments parallel to the stable curve of
the fixed point. This situation happens (in the context of our
homeomorphism of the plane) when the same stable curve of a point
in the unstable curve of the fixed point, is the first one that is
not intersected by the unstable curve of all the points of an arc
of the stable curve of the fixed point. We can approximate $J$ by
curves $J_{n}$ as shown in figure \ref{fig19}. Notice that these
curves $J_{n}$ are similar to those in the previous case.
\begin{figure}[htb]
\psfrag{p1}{$p_{1}$} \psfrag{J1}{$J_{1}$} \psfrag{J}{$J$}
\psfrag{Jn}{$J_{n}$}
\begin{center}
\includegraphics[scale=0.2]{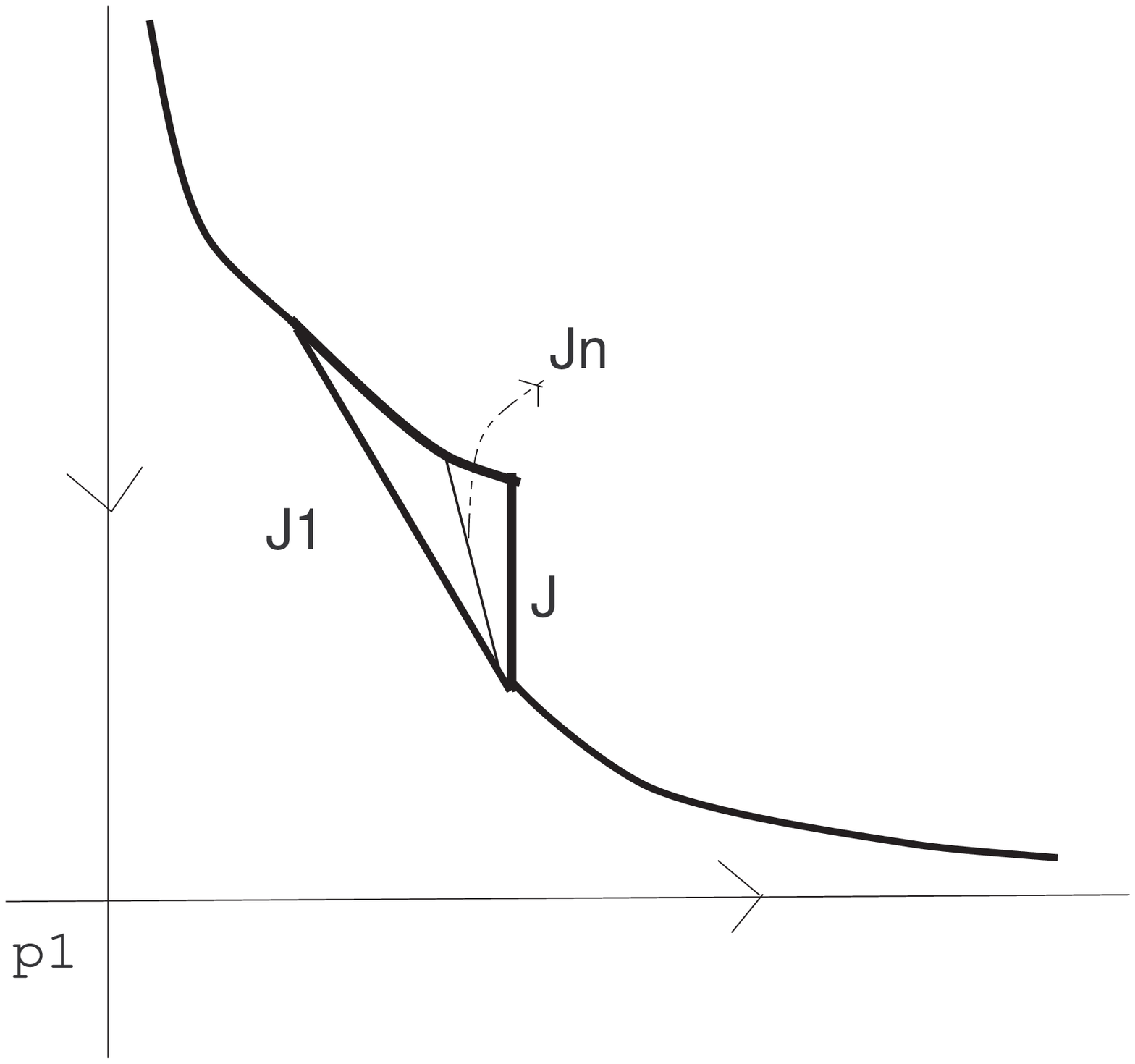}
\caption{\label{fig19}Case 2}
\end{center}
\end{figure}
We can build a conjugation $H$ between the restriction with border
$J$ and the region whose border is $J_{1}$ using the conjugations
$H_{n}$ between the case with border $J_{n}$ and the case with
border $J_{1}$ (to define conjugations $H_{n}$ we use the same
arguments used in case $1$). \item[Case 3] Finally, we will
consider a region whose border admits segments parallel to both
axes, (stable and unstable curves of the fixed point) as shown in
figure \ref{fig20}.
\begin{figure}[htb]
\psfrag{J1}{$J_{1}$} \psfrag{Jn}{$J_{n}$} \psfrag{p1}{$p_{1}$}
\psfrag{J}{$J$}
\begin{center}
\includegraphics[scale=0.2]{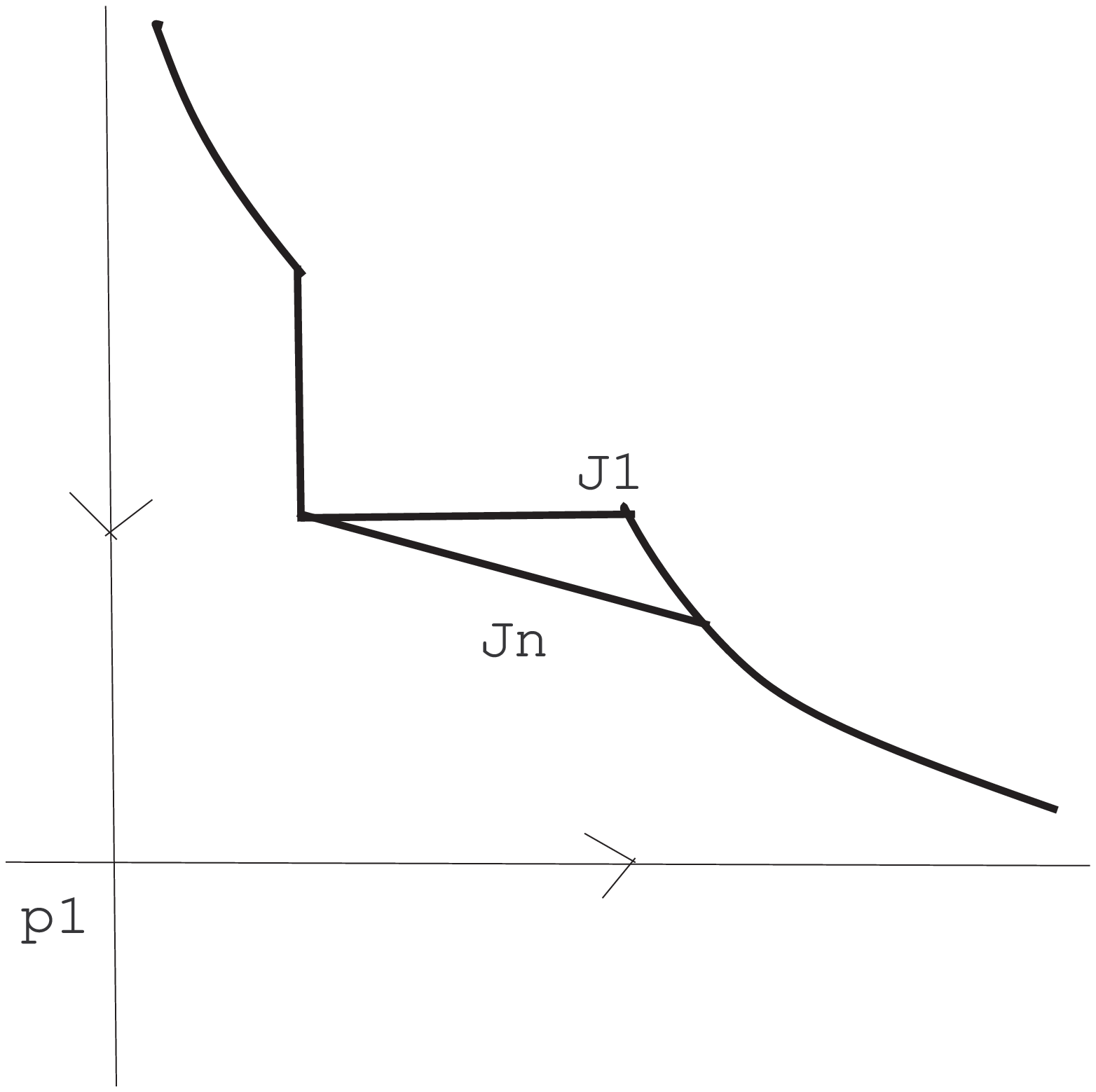}
\caption{\label{fig20}Case 3}
\end{center}
\end{figure}
To prove this case we would use similar arguments to those used in
the previous case.\end{itemize}
\ep

\begin{obs}
Theorem \ref{LOC} refers to a single quadrant. Because of this,
observe that the behavior in each quadrant is independent of the
others and therefore we can obtain different combinations. Notice
that these two classes (referred to a given quadrant) do differ in
the fact that in one class every stable curve intersects every
unstable curve, while in the other one there exist stable curves
that do not intersect some unstable curve.
\end{obs}
Figure \ref{fig59} shows some of these behaviors:
\begin{figure}[htb]
\begin{center}
\includegraphics[scale=0.25]{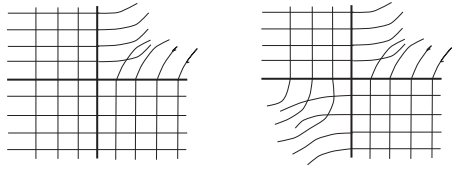}
\caption{\label{fig59}Other behavior.}
\end{center}
\end{figure}
\subsection{Examples.}
\label{EJ} The following two examples show that the classes
determined in Theorem \ref{LOC} are non-empty.

\subsubsection{Example 1.}
We will show an example that admits stable and unstable curves
that do not intersect each other and verifies the conditions shown
in Theorem \ref{LOC}.

{\bf Construction.} Let us consider $f:\R^{2} \rightarrow \R^{2} $
defined by

$$f(x,y)= \left(%
\begin{array}{cc}
  2 & 0 \\
  0 & 1/2 \\
\end{array}%
\right) \left( \begin{array}{c}
  x \\
  y \\
\end{array} \right)    .$$

\begin{figure}[htb]
\psfrag{x0}{$x_{0}$} \psfrag{1x0}{$1/x_{0}$} \psfrag{y0}{$y_{0}$}
\psfrag{1y0}{$1/y_{0}$} \psfrag{Hx0y0}{$H(x_{0},y_{0})$}
\begin{center}
\includegraphics[scale=0.25]{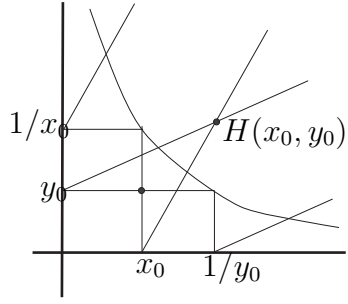}
\caption{\label{fig300}Construction}
\end{center}
\end{figure}

Consider the restriction of $f$ to the region $$\Omega =\{ (x,y):
x\geq 0,\ y\geq 0,\ xy<1\} .$$ We will construct an extension to
the whole quadrant of the linear automorphism restricted to the
region $\Omega $. Each point $(k,1/k)$ on the hyperbola $xy=1$
determines a stable segment parameterized by $(k,u)$, with $0\leq
u\leq 1/k$ and an unstable segment parameterized by $(v,1/k)$,
with $0\leq v\leq k$. Such segments will correspond to semi
straight lines defined by: the semi straight line corresponding to
the stable segment parameterized by $(k,u)$, with $0\leq u\leq
1/k$, defined by $y=1/k(x-k)$ and the semi straight line
corresponding to the unstable segment parameterized by $(v,1/k)$,
with $0\leq v\leq k$, defined by $y=1/k(1+x)$. Let $H:\Omega
\rightarrow \R^{2}$ be such that
$$H(x_{0},y_{0})=\left( \frac{x_{0}(1+y_{0})}{1-x_{0}y_{0}}
,\frac{y_{0}(1+x_{0})}{1-x_{0}y_{0}} \right) ,$$
$H(x_{0},0)=(x_{0},0)$ and $H(0,y_{0})=(0,y_{0})$ (see figure
\ref{fig300}).
Let $\Lambda $ be the first quadrant determined by the stable and
unstable curves of the fixed point. We define $F:\Lambda
\rightarrow \Lambda $ such that $F=HfH^{-1}$. Since $F$ is
conjugated to a restriction of the linear automorphism we can
define, for $F$, the Lyapunov function
$L(p_{1},p_{2})=D(H^{-1}(p_{1}),H^{-1}(p_{2}))$, where
$D=D_{s}+D_{u}$ is the metric Lyapunov function associated to $f$.
The example we built verifies condition {\bf HL} and hypothesis
{\bf HA} since it is conjugated to a restriction of the linear
automorphism (see Theorem \ref{LOC}).

\begin{obs}
The images of the stable and unstable segments of $f$ under $H$
are $L$-stable curves and $L$-unstable curves of $F$. However a
simple calculation shows that this semi-straight lines are not
stable and unstable curves in the sense of the usual metric on the
plane. The following example, will verify that its $L$-stable
curve and $L$-unstable curve are also stable and unstable curves
in the sense of the usual metric.
\end{obs}

\subsubsection{Example 2.}
{\bf Construction.} Let us consider $f:\R^{2} \rightarrow \R^{2} $
defined by
$$f(x,y)= \left(%
\begin{array}{cc}
  2 & 0 \\
  0 & 1/2 \\
\end{array}%
\right) \left( \begin{array}{c}
  x \\
  y \\
\end{array} \right)    .$$

Let us consider the restriction of $f$ to region $$\Omega =\{
(x,y): x\geq 0,\ y\geq 0,\ xy<1\} .$$ We will construct an
extension to the whole quadrant of the linear automorphism
restricted to the region $\Omega $. Each point $(k,1/k)$ on the
hyperbola $xy=1$ determines a stable segment parameterized by
$(k,u)$, with $0\leq u\leq 1/k$ and an unstable segment
parameterized by $(v,1/k)$, with $0\leq v\leq k$. Such segments
will correspond to polygonal lines that will be constructed the
following way:
\begin{itemize}
\item {\bf Case $k\leq 1/2$.} The polygonal line corresponding to
the stable segment parameterized by $(k,u)$, with $0\leq u\leq
1/k$ is defined by: $(k,u)$, when $0\leq u\leq 1/k-1$ and is then
followed by a semi-straight line with slope $1/k$. The polygonal
line corresponding to the unstable segment parameterized by
$(v,1/k)$, with $0\leq v\leq k$ is defined by: $(v,1/k)$, when
$0\leq v\leq k-\frac{k}{1+1/k}$, and is then followed by a
semi-straight line with slope $1/k$ (see figure \ref{fig51}).
\item {\bf Case $k\geq 2$.} For this case we would use similar
arguments to those used in the previous case.

\begin{figure}[htb]
\psfrag{k0}{$k_{0}$} \psfrag{1/2}{$1/2$} \psfrag{2}{$2$}
\psfrag{k}{$k$} \psfrag{y=1/k}{$y=1/k$} \psfrag{1}{$1$}
\begin{center}
\includegraphics[scale=0.3]{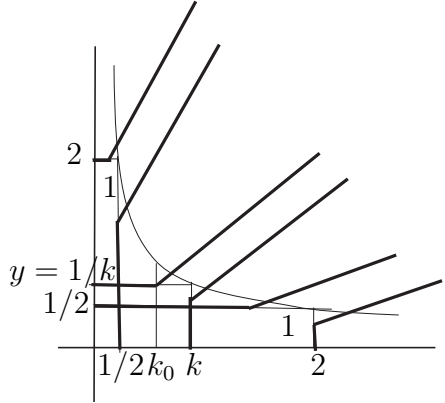}
\caption{\label{fig51}Construction}
\end{center}
\end{figure}

\item {\bf Case $1/2\leq k\leq 2$.} Note that the coordinates of
the point where the polygonal line corresponding to the stable
segment parameterized by $(1/2,u)$, with $0\leq u\leq 2$, breaks
is $(1/2,2-1)$ and it is point $(2,1/2-1/6)$ for the case
corresponding to the stable segment $(2,v)$ with $0\leq v\leq
1/2$. Let $f:[1/2,2]\rightarrow \R$ be the linear function such
that $f(1/2)=1$ y $f(2)=1/6$. The polygonal line corresponding to
the stable segment parameterized by $(k,u)$ is defined by:
$(k,u)$, when $0\leq u\leq 1/k-f(k)$, and is then followed by a
semi-straight line with slope $1/k$. Let $k_{0}<k$ be such that
the length of the vertical segment included in the line $x=k_{0}$
between $y=1/k$ and the hyperbola $xy=1$ is $f(k_{0})$. The
polygonal line corresponding to the unstable segment parameterized
by $(v,1/k)$, with $0\leq v\leq k$ is defined by: $(v,1/k)$ for
$0\leq u\leq k_{0}$, and is then followed by a semi straight-line
with slope $1/k$ (see figure \ref{fig51}).
\end{itemize}
This way, we can define a homeomorphism $H$ that takes the
invariant region $\Omega $ in the first quadrant $\Lambda $,
sending the intersection of a stable segment with an unstable
segment into the intersection of the corresponding polygonal. Let
$F:\Lambda \rightarrow \Lambda $ be such that $F=HfH^{-1}$. Since
$F$ is conjugated to a restriction of the linear $f$, we can
define a Lyapunov metric function for $F$,
$L(p_{1},p_{2})=D(H^{-1}(p_{1}),H^{-1}(p_{2}))$, where
$D=D_{s}+D_{u}$ is the Lyapunov function associated to $f$. The
constructed example verifies condition {\bf HL} and hypothesis
{\bf HA} since it is conjugated to a restriction of a linear
automorphism (preserving stable and unstable curves) (see Theorem
\ref{LOC}).

\begin{obs}
Polygonal lines constructed in this example are stable (unstable)
curves of $F$ not only on a Lyapunov function $L$ sense, but also
in the usual metric sense.
\end{obs}

\begin{figure}[htb]
\psfrag{A}{$A$} \psfrag{B}{$B$} \psfrag{H-1p}{$H^{-1}(p)$}
\psfrag{H-1q}{$H^{-1}(q)$} \psfrag{p}{$p$} \psfrag{q}{$q$}
\psfrag{1}{$1$}
\begin{center}
\includegraphics[scale=0.25]{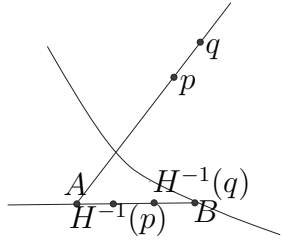}
\caption{\label{fig53}Stable and unstable}
\end{center}
\end{figure}

\bp Let us take as an example two points $p,q$ that are in the
same unstable polygonal line, see figure \ref{fig53}. Because of
the construction built before, the length of the unstable segment
$AB$ tends to zero for the past with the same order that
$\frac{1}{y^{2}} $ when $y$ tends to infinity. While the length of
segment $H^{-1}(p),H^{-1}(q)$ does it with order $1/y$.

Then, iterating for the past, points
$f^{-n}(H^{-1}(p)),f^{-n}(H^{-1}(q))$ get inside the zone where
$H$ is the identity and therefore the usual distance between
$F^{-n}(p),F^{-n}(q)$ tends to zero when $n$ tends to
infinity.\ep

\end{document}